\documentclass[10pt]{article}
\usepackage{amstext}
\usepackage{amsfonts}
\usepackage{amssymb}
\usepackage{amsbsy}
\usepackage{amsmath}
\usepackage{latexsym}
\usepackage{xy}
\usepackage{hhline}
\xyoption{all}
\newcommand{\vtx}[1]{*+[o][F-]{\scriptscriptstyle #1}} 

\newcounter{num}[section] %

\newenvironment{theo}
{\refstepcounter{num}%
\bigskip\noindent{\bf Theorem~\arabic{section}.\arabic{num}. }\it}

\newenvironment{cor}
{\refstepcounter{num}%
\bigskip\noindent{\bf Corollary~\arabic{section}.\arabic{num}. }\it}

\newenvironment{lemma}
{\refstepcounter{num}%
\bigskip\noindent{\bf Lemma~\arabic{section}.\arabic{num}. }\it}

\newenvironment{example}
{\refstepcounter{num}%
\bigskip\noindent{\bf Example~\arabic{section}.\arabic{num}.}}

\newenvironment{remark}
{\refstepcounter{num}%
\bigskip\noindent{\bf Remark~\arabic{section}.\arabic{num}.}}

\newcommand{\defin}{\refstepcounter{num}%
\bigskip\noindent{\bf Definition~\arabic{section}.\arabic{num}.}}

\newcommand{\Ref}[1]{(\ref{#1})}

\newcounter{thepic}

\newenvironment{proof}{\medskip\noindent{\it Proof. }}
{$\Box$ \bigskip}
\newenvironment{proof_of}[1]{\medskip\noindent{\it Proof #1}}
{$\Box$ \bigskip}
\newenvironment{eq}{\begin{equation}}{\end{equation}}

\newcommand{\si}{\sigma}
\newcommand{\al}{\alpha}
\newcommand{\be}{\beta}
\newcommand{\ga}{\gamma}
\newcommand{\la}{\lambda}
\newcommand{\de}{\delta}
\newcommand{\De}{\Delta}
\newcommand{\La}{\Lambda}
\newcommand{\Ga}{\Gamma}

\newcommand{\ov}[1]{\overline{#1}}

\newcommand{\tr}{\mathop{\rm tr}}

\newcommand{\mdeg}{\mathop{\rm mdeg}}

\newcommand{\Char}{\mathop{\rm char}}
\newcommand{\sign}{\mathop{\rm{sgn }}}

\newcommand{\Hom}{{\mathop{\rm{Hom }}}}
\newcommand{\supp}{{\mathop{\rm{supp }}}}

\newcommand{\algA}{\mathcal{A}}    

\newcommand{\FF}{{\mathbb{F}}}   
\newcommand{\NN}{{\mathbb{N}}}

\newcommand{\Q}{\mathcal{Q}}    
\newcommand{\n}{\boldsymbol{n}} 
\newcommand{\Ver}[1]{\mathop{{\rm ver}(#1)}} 
\newcommand{\Arr}[1]{\mathop{{\rm arr}(#1)}} 
\newcommand{\Inv}{I}       
\newcommand{\path}[1]{\mathop{{\rm path}(#1)}} 

\newcommand{\SetII}{\mathcal{S}_2}                 
\newcommand{\SetI}{\mathcal{S}_1}                 
\newcommand{\loopR}[3]{%
\begin{picture}(20,0)(#1,#2)
\put(-2,1){\llap{$\scriptstyle #3$}} \put(10,3){\circle{20}} \put(20,6){\vector(1,-4){1}}
\end{picture}}
\newcommand{\loopL}[3]{%
\begin{picture}(20,0)(#1,#2)
\put(22,1){$\scriptstyle #3$} \put(10,3){\circle{20}} \put(0,6){\vector(-1,-4){1}}
\end{picture}}

\newcommand{\rectangle}[2]{
\begin{picture}(0,0)
\put(-#1,-#2){\line(1,0){#1}}\put(0,-#2){\line(1,0){#1}}
\put(-#1,#2){\line(1,0){#1}}\put(0,#2){\line(1,0){#1}}
\put(-#1,-#2){\line(0,1){#2}}\put(-#1,0){\line(0,1){#2}}
\put(#1,-#2){\line(0,1){#2}}\put(#1,0){\line(0,1){#2}}
\end{picture}}

\begin{document}
\title{Minimal generating set for semi-invariants of quivers of dimension two.}
\author{
A.A. Lopatin\thanks{Supported by DFG} \\
{\small\it Institute of Mathematics, SBRAS, Pevtsova street, 13, Omsk 644099, Russia} \\
{\small\it artem\underline{ }lopatin@yahoo.com} \\
}
\date{} 
\maketitle

\begin{abstract}
A minimal (by inclusion) generating set for the algebra of semi-invariants of a quiver of dimension $(2,\ldots,2)$ is established over an infinite field of arbitrary
characteristic. The mentioned generating set consists of the determinants of generic matrices and the traces of tree paths of pairwise different multidegrees, where in the case of characteristic different from two we take only admissible paths.  As a consequence, we describe relations modulo decomposable semi-invariants.
\end{abstract}




\section{Introduction}\label{section_intro}
We work over an infinite field $\FF$ of arbitrary characteristic $\Char(\FF)$. All vector
spaces, algebras, and modules are over $\FF$ and all algebras are
associative unless otherwise stated.

A {\it quiver} $\Q=(\Q_0,\Q_1)$ is a finite oriented graph, where $\Q_0$ stands for the set of vertices and $\Q_1$ stands for the set of arrows. For an arrow $a$ denote by $a'$
its head and denote by $a''$ its tail.  The notion of quiver was introduced by Gabriel in~\cite{Gabriel_1972} and it was applied to describe different problems of the linear algebra. The importance of this notion from point of
view of the representation theory is due to the following fact. Let $\algA$ be a finite dimensional basic algebra over an 
algebraically closed field. Then the category of finite dimensional modules over $\algA$
is a full subcategory of the category of representations of some
quiver (for example, see Chapter~3 from~\cite{Kirichenko}). 

Given a {\it dimension vector}
$\n=(\n_v\,|\,v\in\Q_0)$, we assign an $\n_v$-dimensional vector space
$V_v$ to $v\in \Q_0$. We identify $V_v$ with the space of column
vectors $\FF^{\n_v}$. Fix the {\it standard} basis
$e(v,1),\ldots,e(v,\n_v)$ for $\FF^{\n_v}$, where $e(v,i)$ is a
column vector whose $i^{\rm th}$ entry is $1$ and the rest of entries
are zero. A {\it representation} of $\Q$ of dimension vector
$\n$ is a collection of matrices %
$$h=(h_a)_{a\in \Q_1}\in %
H=H(\Q,\n)=\bigoplus_{a\in \Q_1} \FF^{\n_{a'}\times \n_{a''}} \simeq%
\bigoplus_{a\in \Q_1} \Hom_\FF(V_{a''},V_{a'}),$$ %
where $\FF^{n_1\times n_2}$ stands for the linear space of
$n_1\times n_2$ matrices over $\FF$ and the isomorphism is given by
the choice of bases. We will refer to $H$ as the {\it space of
representations} of $\Q$ of dimension vector $\n$.
The group %
$$GL(\n)=\prod_{v\in \Q_0} GL(\n_v)$$ %
acts on $H$ as change of the bases for $V_v$ ($v\in \Q_0$). In
other words, $GL(\n_v)$ acts on $V_v$ by left multiplication, and
this action induces the action of $G$ on $H$ by
$$g\cdot h=(g_{a'}h_a g_{a''}^{-1})_{a\in \Q_1},$$ %
where $g=(g_a)_{a\in \Q_1}\in GL(\n)$ and $h=(h_a)_{a\in \Q_1}\in H$.

The coordinate ring of the affine variety $H$ is the polynomial ring  %
$$\FF[H]=\FF[x_{ij}^a\,|\,a\in \Q_1,\,1\leq i\leq \n_{a'},1\leq j\leq\n_{a''}],$$ %
where $x_{ij}^a$ is the coordinate function on $H$ that
takes a representation $h\in H$ to the $(i,j)^{\rm th}$ entry of a
matrix $h_{a}$. Denote by
$$X_a=\left(\begin{array}{ccc}
x_{1,1}^a & \cdots & x_{1,n_{a''}}^a\\
\vdots & & \vdots \\
x_{n_{a'},1}^a & \cdots & x_{n_{a'},n_{a''}}^a\\
\end{array}
\right)
$$%
the $\n_{a'}\times \n_{a''}$ {\it
generic} matrix. The action of $GL(\n)$ on $H$ induces the action on $\FF[H]$ as follows:
$(g\cdot f)(h)=f(g^{-1}\cdot h)$ for all $g\in GL(\n)$, $f\in \FF[H]$,
$h\in H$. In other words, %
\begin{eq}\label{eq_action}
g\cdot x^a_{ij}=(i,j)^{\rm th}\text{ entry of }g_{a'}^{-1}X_a g_{a''}.
\end{eq} %
The algebra of {\it invariants} is
$$I(\Q,\n)=\FF[H]^{GL(\n)}=\{f\in \FF[H]\,|\,g\cdot f=f\;{\rm for\; all}\;g\in GL(\n)\}.$$
Similarly, for the group
$$SL(\n)=\prod_{v\in \Q_0} SL(n_v) < GL(\n)$$
we define the algebra of {\it semi-invariants}
$$SI(\Q,\n)=\FF[H]^{SL(\n)}.$$

To describe generators for $I(\Q,\n)$ we use the following notions.
Denote by $\si_t(X)$ the $t^{\rm th}$ coefficient
in the characteristic polynomial of an $n\times n$ matrix $X$,
i.e., %
$$\det(\la E+X)=\la^n+\si_1(X)\la^{n-1}+\cdots+\si_n(X).
$$%
In particular, $\si_1(X)=\tr(X)$ and $\si_n(X)=\det(X)$. We say that $a=a_1\cdots a_s$ is a {\it path} in $\Q$ (where $a_1,\ldots,a_s\in
\Q_1$), if $a_1''=a_2',\ldots,a_{s-1}''=a_s'$, i.e.,
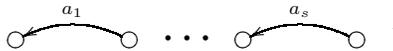
$$\vcenter{
\xymatrix@C=1.3cm@R=1.3cm{ %
\vtx{}\ar@/^/@{<-}[r]^{a_1} & %
\vtx{} & %
\vtx{}\ar@/^/@{<-}[r]^{a_s} & %
\vtx{}\\
}} \quad.
\begin{picture}(0,0) 
\put(-80,-3){
\put(0,0){\circle*{2}} %
\put(-7,0){\circle*{2}} %
\put(7,0){\circle*{2}} %
} %
\end{picture}$$
The head of the path
$a$ is $a'=a_1'$ and the tail of $a$ is $a''=a_s''$. If $a_1'=a_s''=v$, then the path $a$ is called {\it closed} in the vertex $v$. We write $X_{a}$ for $X_{a_1}\cdots X_{a_s}$. 
Denote the degree of $a$ by $\deg(a)=s$ and the degree of $a$ in an arrow $x$ by $\deg_{x}(a)$. 

In~\cite{Donkin_1994} Donkin proved that $\FF$-algebra $I(\Q,\n)$ is generated by
$\si_t(X_{a})$ for all closed paths $a$ in $\Q$ and $1\leq t\leq n_{a'}$. In characteristic zero case 
generators for $\Inv(\Q,\n)$ had earlier been described by Le Bruyn and Procesi
in~\cite{Le_Bruyn_Procesi_1990}. Relations between generators were described by Zubkov in~\cite{Zubkov_Fund_Math_2001}. In characteristic zero case 
this result had earlier been obtained by Domokos in~\cite{Domokos_1998}.

Generators for $SI(\Q,\n)$ were described by
Domokos and Zubkov in~\cite{DZ01}  and,
independently, by Derksen and Weyman in~\cite{DW00},~\cite{DW_LR_02}. Simultaneously, similar result in the
case of characteristic zero was obtained by Schofield and van den Bergh
in~\cite{Schofield_van_den_Bergh_2001}. These results were generalized for semi-invariants of mixed representations of quivers by
the author and Zubkov in~\cite{LZ1} and for semi-invariants of supermixed representations by the author in~\cite{Lopatin_so_inv}.

In this paper we assume that $\n=(2,\ldots,2)$ unless otherwise stated and write $I(\Q)$,  $SI(\Q)$ for $I(\Q,\n)$ and $SI(\Q,\n)$, respectively.

\defin{} Define the quiver $\Q^{\ast}$ as follows: $\Q^{\ast}_0=\Q_0$ and $\Q^{\ast}_1=\Q_1 \sqcup \{a^{\ast}\,|\,a\in \Q_1\}$, where $(a^{\ast})'=a''$ and $(a^{\ast})''=a'$.
We set $X_{a^{\ast}}=X_{a}^{\ast}=-J_2 X^{T}_{a}J_2$ for all $a\in\Q_1$, where $J_2=\left(
\begin{array}{cc}
0& 1 \\
-1& 0\\
\end{array}
\right)$ %
is the matrix of the skew-symmetric bilinear form on $\FF^2$.  
\bigskip

Let $g\in SL(\n)$ and $a\in\Q_1$. By~\Ref{eq_action}, $\det(X_a)$ is a semi-invariant.  For short, we write $g\cdot X_{a}$ for the matrix whose $(i,j)^{\rm th}$ entry is $g\cdot x_{ij}^a$. Defining $g\cdot X_{a^{\ast}}$ similarly, we have 
\begin{eq}\label{eq_action_ast}
g\cdot X_{a^{\ast}}=-J_2 (g\cdot X_a)^T J_2 = g_{a''}^{-1} X_{a^{\ast}} g_{a'},
\end{eq}%
where we use the equality $AJ_2A^{T}=J_2$ for an $A\in SL(2)$. It is not difficult to see that~\Ref{eq_action} and~\Ref{eq_action_ast} imply that $\tr(X_b)$ is a semi-invariant for any closed path $b$ in $\Q^{\ast}$.  

Applying the general description of semi-invariants from~\cite{DZ01}, Fedotov has recently showed that in characteristic zero case the above mentioned semi-invariants actually generate the algebra $SI(\Q)$  (see~\cite{Fedotov}). In this paper we present an independent proof, which also covers the case of arbitrary characteristic. 

\begin{theo}\label{theo_Fedotov} The algebra $SI(\Q)$ is generated by 
\begin{enumerate}
\item[$\bullet$] $\det(X_a)$ and $\tr(X_b)$, if $\Char{\FF}=2$; 

\item[$\bullet$] $\tr(X_b)$, otherwise,
\end{enumerate}
where $a$ ranges over $\Q_1$ and $b$ ranges over all closed paths $b$ in $\Q^{\ast}$. Moreover, we can assume that $\deg_{x}(b)\leq 1$ for all $x\in \Q_1^{\ast}$. 
\end{theo}
\bigskip

\noindent{}Our main result is the explicit description of a minimal (by inclusion) generating set for the $\FF$-algebra $SI(\Q)$: 
\begin{enumerate}
\item[$\bullet$] if $\Char{\FF}=2$, then see Theorem~\ref{theo_mgs_2}; note that in the case of arbitrary characteristic Theorem~\ref{theo_mgs_2} yields the generating set for $SI(\Q)$, which is smaller than the generating set from Theorem~\ref{theo_Fedotov} and is not as complicated as the generating set from Theorem~\ref{theo_mgs_not2};    

\item[$\bullet$] if $\Char{\FF}\neq2$, then see Theorem~\ref{theo_mgs_not2}.
\end{enumerate}

Let $\NN=\{0,1,2,\ldots\}$ and $\n$ be arbitrary. The algebra $\FF[H]$ has the natural $\NN$-grading by degrees and $\NN^{\#\Q_1}$-grading by multidegrees defined as follows. For a monomial $f\in \FF[H]$ we set 
$$\deg(f)=\sum_{a\in\Q_1} \deg_a(f) \;\text{ and }\; \mdeg(f)=(\deg_a(f)\,|\,a\in\Q_1),$$  
where $\deg_a(f)=\sum_{ij} \deg_{x_{ij}^a}(f)$.   The algebras $I(\Q,\n)$ and  $SI(\Q)$ have also gradings induced by the mentioned ones. To describe these gradings for the given generators of algebras, we introduce 
the multidegree of a path $b$ in $\Q^{\ast}$ by $\mdeg(b)=(d_a\,|\,a\in \Q_1)$, where $d_a=\deg_a(b)+\deg_{a^{\ast}}(b)$. Since  
$$\deg_a(\si_t(X_b))=t(\deg_a(b)+\deg_{a^{\ast}}(b))$$
for an $a\in\Q_1$ and $t>0$, we have
$$\deg(\si_t(X_{b}))=t \deg(b)\;\text{ and }\; \mdeg(\si_t(X_b))=t\mdeg(b).$$
Note that minimal generating sets from Theorems~\ref{theo_mgs_2} and~\ref{theo_mgs_not2} are $\NN^{\#\Q_1}$-homogeneous. 

Given an $\NN$-graded algebra $A$, denote by $A^{+}$ the subalgebra generated by elements of $A$ of positive degree. It is easy to see that an $\NN$-homogeneous set $\{a_i\} \subseteq A$ is a minimal set of generators if and only if $\{\ov{a_i}\}$ is a basis of $\ov{A}={A}/{(A^{+})^2}$. We say that an element $a\in A$ is {\it decomposable} and write $a\equiv0$ if it belongs to the ideal $(A^{+})^2$. In other words, a decomposable element is equal to a polynomial in elements of strictly lower degree.

As a consequence of Theorems~\ref{theo_mgs_2} and~\ref{theo_mgs_not2}, we obtain the following results:
\begin{enumerate}
\item[$\bullet$] Up to isomorphism, $SI(\Q)$ does not depend on the orientation of arrows of $\Q$ (see Corollary~\ref{cor_orientation}).

\item[$\bullet$] Up to multiplication on elements of $\FF$, a minimal $\NN^{\#\Q_1}$-homogeneous generating set for $SI(\Q)$ is unique modulo indecomposable semi-invariants (see Corollary~\ref{cor_MGS}).  

\item[$\bullet$] Relations between generators for $SI(\Q)$ are described modulo decomposable semi-invariants (see Corollary~\ref{cor_relations}). Note that the ideal of relations between generators for $I(\Q)$ is known in contrast to $SI(\Q)$. Nevertheless,  the only known result concerning a minimal generating set for $I(\Q)$ is an upper bound on degrees of indecomposable invariants (see~\cite{Lopatin_2222_I} and~\cite{Lopatin_2222_II}) and a minimal generating set for $I(\Q)$ is still not known. 

\item[$\bullet$] As an example, we consider a partial case of so-called tree-like quivers in Section~\ref{section_example}. In particular, we prove that if $\Q$ is a tree, then $SI(\Q)$ is a polynomial algebra. Considering a quiver with two vertices, we compare the generating set from Theorem~\ref{theo_Fedotov} with the minimal generating set from Theorems~\ref{theo_mgs_2} and~\ref{theo_mgs_not2} (see Remark~\ref{remark_compare}).
\end{enumerate}

The paper is organized as follows. Section~\ref{section_results} contains formulations of the main results. Using the notion of tableau with substitution introduced in~\cite{Lopatin_bplp}, we prove Theorem~\ref{theo_Fedotov} in Section~\ref{section_generating_set}. Key lemmas are proven in Section~\ref{section_relations}. The proof of Theorem~\ref{theo_mgs_2} is given at the end of Section~\ref{section_case2} and Theorem~\ref{theo_mgs_not2} is proven at the end of Section~\ref{section_case_not2}. Some applications and examples are considered in Sections~\ref{section_corollaries} and~\ref{section_example}, respectively.

\section{Results}\label{section_results}

\subsection{The case of $\Char{\FF}=2$}\label{section_results_char2}

We start this section with some definitions. Let $a=a_1\cdots a_s$ be a path in $\Q^{\ast}$, where $a_1,\ldots,a_s\in\Q_1^{\ast}$. Denote $\Ver{a}=\{a_1',\ldots,a_s',a_s''\}$ and $\Arr{a}=\{a_1,\ldots,a_s\}$. We define $\Ver{\Omega}$ and $\Arr{\Omega}$ for a set $\Omega\subset\Q_1^{\ast}$ similarly. We write $(a^{\ast})^{\ast}$ for $a$ and $a^{\ast}$ for the path $a_s^{\ast} \cdots a_1^{\ast}$ in $\Q^{\ast}$. Note that $X_{a^{\ast}}=X_{a}^{\ast}$ and $X_{(a^{\ast})^{\ast}}=(X_{a}^{\ast})^{\ast}=X_{a}$ 
for any path $a$. Denote by $\supp(a)$ the {\it support} of $a$, i.e., $\supp(a)$ is a quiver with $\supp(a)_0=\Ver{a}$ and $\supp(a)_1=\Arr{a}$.

\defin{}\label{def_miltilinear_path} A {\it multilinear path} in $\Q^{\ast}$ is a closed path $a$ in $\Q^{\ast}$ such that $\deg_x(a)\leq1$ for all $x\in\Q_1^{\ast}$.

\defin{}\label{def_tree_path} A {\it tree} path $a$ is a multilinear path in $\Q^{\ast}$ such that if $\deg_x(a)=\deg_{x^{\ast}}(a)=1$ for an $x\in\Q_1$, then  $$\Arr{a}=\Gamma_1\sqcup \Lambda_1\sqcup\{x,x^{\ast}\}$$
for some quivers $\Ga,\La$ with non-empty $\Ga_1$, $\La_1$, and $\Ga_0\cap\La_0=\emptyset$, i.e.,
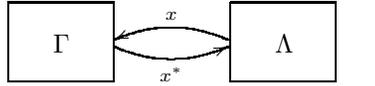
$$\xymatrix@C=1.5cm@R=1.5cm{ %
\ar@/^/@{<-}[r]^{x}  \ar@/_/@{->}[r]_{x^{\ast}} &} %
\begin{picture}(0,0)
\put(-69,3){\rectangle{20}{15}\put(-3,-4){$\Ga$}}%
\put(15,3){\rectangle{20}{15}\put(-3,-4){$\La$}}%
\put(50,0){.}
\end{picture}
$$

Let $a$ be a tree path in $\Q^{\ast}$. Then %
$$\Ver{a}=\bigsqcup_{j\in J}\, \Gamma^{(j)}_0 \;\text{ and }\; \Arr{a}=\bigsqcup_{i\in I}\, \{x_i,x_i^{\ast}\} \sqcup \bigsqcup_{j\in J}\, \Gamma^{(j)}_1,$$
where $x_i\in \Q_1$ and $\Gamma^{(j)}$ is a quiver with non-empty set of vertices such that  for any $j$ there is no $x\in \Q_1$ with $\{x,x^{\ast}\}\subset \Gamma^{(j)}_1$. Moreover, consider a graph with vertices $1,\ldots,\#J$ and edges $1,\ldots,\#I$, where an edge $i$ connects vertices $j_1$ and $j_2$ if and only if  $x_i'\in \Gamma^{(j_1)}_0$ and $x_i''\in \Gamma^{(j_2)}_0$. Hence this graph is a tree and it is called the {\it tree} of path $a$. 
Quivers $\Gamma^{(j)}$ are called {\it blocks} of $a$. Note that some blocks can have empty sets of arrows. But if a block corresponds to a leaf of tree of path $a$, then the set of arrows of this block is not empty. 

We denote by $\SetII$ a maximal (by inclusion) subset of tree paths such that elements of $\SetII$ have pairwise different multidegrees.

\begin{theo}\label{theo_mgs_2}
The algebra $SI(\Q)$ is generated by $\{\det(X_a), \tr(X_b)\,|\, a\in\Q_1,\;b\in \SetII\}$ over $\FF$. 

Moreover, if $\Char{\FF}=2$, then the given set is a minimal generating set for $SI(\Q)$.
\end{theo}

\subsection{The case of $\Char{\FF}\neq2$}\label{section_results_char0}

Let $\Char(\FF)\neq2$.  Since $\tr(X_{a_1}\cdots X_{a_4})\equiv0$  
for all closed paths $a_1,\ldots,a_4$ in $\Q^{\ast}$ with $a_1'=\cdots=a_4'$ (see Lemma~\ref{lemma_two_rel} below), we have to remove some elements from the generating set given in Theorem~\ref{theo_mgs_2} to obtain a minimal generating set.  
To perform this operation we introduce the following notions. 

We endow the set of closed paths in $\Q^{\ast}$ with the equivalence $\sim$ as follows:
\begin{enumerate}
\item[$\bullet$] $a\sim a^{\ast}$,

\item[$\bullet$] $xy\sim yx$,
\end{enumerate}
where $a$ and $xy$ are closed paths in $\Q^{\ast}$. As an example, if $x y^{\ast} z$ is a closed path in $\Q^{\ast}$, then 
$$x y^{\ast} z \sim x^{\ast} z^{\ast} y.$$

\defin{}\label{def_decomposition} Assume that $a$ is a closed path in $\Q^{\ast}$. We say that $\{b_1,\ldots,b_s\}$ is a {\it decomposition} of $a$ into primitive closed paths, if 
\begin{enumerate}
\item[$\bullet$] $b_1,\ldots,b_s$ are primitive (i.e., without self-intersections) closed paths in $\Q^{\ast}$ satisfying $\Arr{b_i}\cap \Arr{b_j}=\emptyset$ for $i\neq j$;

\item[$\bullet$] $\Arr{a}=\Arr{b_1} \sqcup\cdots \sqcup\Arr{b_s}$.
\end{enumerate}

\begin{remark} A decomposition of a closed path into primitive closed paths is not unique (see Example~\ref{ex_decomp} below).
\end{remark}

\example{}\label{ex_decomp} 
Let $\Q$ be the following quiver:
$$\xymatrix@C=1.3cm@R=1.3cm{ %
\vtx{} \ar@/^/@{->}[r]^{a_3} \ar@/_/@{<-}[d]_{a_4} & %
\vtx{} \ar@/_/@{<-}[d]_{b_4} \ar@/^/@{->}[d]^{a_2} \ar@/^/@{->}[r]^{b_3} & %
\vtx{} \ar@/_/@{<-}[d]_{c_4} \ar@/^/@{->}[d]^{b_2} \ar@/^/@{->}[r]^{c_3} & %
\vtx{} \ar@/^/@{->}[d]^{c_2} \ar@/_/@{<-}[d]_{d_4} \ar@/^/@{->}[r]^{d_3} &
\vtx{} \ar@/^/@{->}[d]^{d_2} \\
\vtx{} \ar@/_/@{<-}[r]_{a_1} & \vtx{} \ar@/_/@{<-}[r]_{b_1} & %
\vtx{} \ar@/_/@{<-}[r]_{c_1} & \vtx{} \ar@/_/@{<-}[r]_{d_1} & \vtx{}\\
}$$%
\noindent{}Then the closed path $h=a_1 a_2 b_4 b_1 b_2 c_4 c_1 c_2 d_4 d_1 d_2 d_3 c_3 b_3 a_3 a_4$ in $\Q$ is multilinear and it has the following decompositions into primitive closed paths: 
\begin{enumerate}
\item[a)] $\{a,b,c,d\}$ for $a=a_1 b_1 c_1 d_1 d_2 d_3 c_3 b_3 a_3 a_4$, $b=a_2 b_4$, $c=b_2 c_4$, $d=c_2 d_4$;

\item[b)] $\{a,b,c,d\}$ for $a=a_1 a_2 a_3 a_4$, $b=b_1 b_2 b_3 b_4$, $c=c_1 c_2 c_3 c_4$,  $d=d_1 d_2 d_3 d_4$.
\end{enumerate}

%

\begin{lemma}\label{lemma_decomposit} Every multilinear path in $\Q^{\ast}$ has some decomposition into primitive closed paths.
\end{lemma}
\begin{proof} Let $a=a_1\cdots a_r$ be a closed multilinear path in $\Q^{\ast}$, where $a_1,\ldots,a_r\in\Q_1^{\ast}$. We prove the lemma by induction on $r\geq 1$. 

If $r=1$, then $\{a_1\}$ is the required decomposition.

Let $r>1$. Then one of the following possibilities holds.
\begin{enumerate}
\item[a)] $a$ is primitive. Then $\{a\}$ is the required decomposition.

\item[b)] $a=xby$, where $b$ is a primitive closed path in $\Q^{\ast}$, $y$ is a path, and  $x=x_1\cdots x_l$ for $x_1,\ldots,x_l  \in\Q_1^{\ast}$ with pairwise different vertices $x_1',\ldots, x_l',x_l''$. By induction hypothesis, there exists a decomposition $\{b_1,\ldots,b_s\}$ of $xy$ into primitive closed paths. Since $a$ is multilinear, $\{b,b_1,\ldots,b_s\}$ is the required decomposition.    

\item[c)] $a=by$, where $b$ is a primitive closed path in $\Q^{\ast}$ and $y$ is a path. This case is similar to case~b).
\end{enumerate}
\end{proof} 

\defin{}\label{def_diagram} A {\it diagram} $D$ is a finite graph without loops and multiple edges such that its edges are marked with positive integers. 

\defin{}\label{def_type} 
Assume that $a$ is a multilinear path in $\Q^{\ast}$, $\{b_1,\ldots,b_s\}$ is some  decomposition of $a$ into primitive closed paths. Then the following diagram $D$ is called the {\it type of $a$ with respect to} $\{b_1,\ldots,b_s\}$:
\begin{enumerate}
\item[a)] $D_0=\{v_1,\ldots,v_s\}$.

\item[b)] If $b_i$ and $b_j$ do not intersect, i.e., $\#\Ver{b_i}\cap\Ver{b_j}=0$, then there is no edge in $D$ that goes from $v_i$ to $v_j$ ($i\neq j$). 

\item[c)] If $b_i$ and $b_j$ intersect at $t$ different vertices, i.e., $\#\Ver{b_i}\cap\Ver{b_j}=t$, then there is an edge in $D$ that goes from $v_i$ to $v_j$ ($i\neq j$) and this edge is marked with $t$.
\end{enumerate}

\example{}\label{ex_type} Let $h$ be the closed path from Example~\ref{ex_decomp} and $\{a,b,c,d\}$ be the decomposition of $h$ from part~a) (part~b), respectively) of Example~\ref{ex_decomp}. Then the type of $h$ with respect to $\{a,b,c,d\}$ is the following diagram $B$ ($D$, respectively):
$$B:\quad
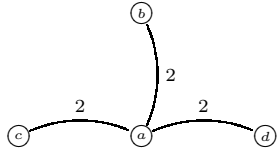
\xymatrix@C=1.3cm@R=1.3cm{ %
& \vtx{b} \ar@/^/@{-}[d]^{2} & \\
\vtx{c}   \ar@/^/@{-}[r]^{2} &\vtx{a}\ar@/^/@{-}[r]^{2} & \vtx{d} \\
}
\qquad\qquad\qquad D:\quad
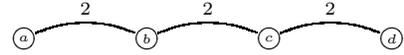
\xymatrix@C=1.3cm@R=2cm{ %
\vtx{a} \ar@/^/@{-}[r]^{2} & \vtx{b} \ar@/^/@{-}[r]^{2} & %
\vtx{c} \ar@/^/@{-}[r]^{2} & \vtx{d}}%
$$%
Here vertices of $B$ and $D$ are denoted by the corresponding closed primitive paths from the given decomposition. 
 
\defin{}\label{def_diagram_admissible} A diagram $D$ is called {\it admissible} if  
\begin{enumerate}
\item[$\bullet$] its edges are marked only with $1,2$; 

\item[$\bullet$] if $a$ is a closed primitive path in $D$, then $\deg{a}=3$ and all edges of $a$ are marked with $1$; in this case we say that $a$ is a {\it triangle}.  
\end{enumerate}

\begin{remark}\label{remark_shrink} Note that every two different triangles of an admissible diagram $D$ do not have a common edge. Moreover, if we shrink all triangles, then $D$ turns into a tree, where it is said that we {\it shrink} a triangle, if we remove its edges and add a new vertex $u$ and new edges that connect $u$ with vertices of this triangle, i.e.,
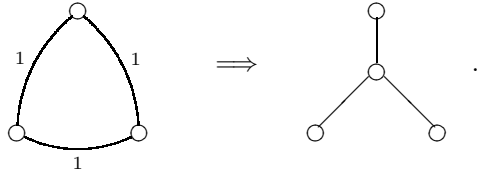
$$\xymatrix@C=0.6cm@R=0.6cm{ %
&\vtx{} \ar@/^/@{-}[rdd]^{1} \ar@/_/@{-}[ldd]_{1} &  \\
&&\\
\vtx{} \ar@/_/@{-}[rr]_{1} && \vtx{} \\
}  \qquad \qquad \qquad 
\xymatrix@C=0.6cm@R=0.6cm{ %
&\vtx{} &  \\
&\vtx{} &\\
\vtx{}  && \vtx{} \\
}   
\begin{picture}(0,0)
\put(-87,-20){$\Longrightarrow$}
\put(10,-20){.}
\put(-26,-17){\line(0,1){17}}
\put(-29,-22){\line(-1,-1){19}}
\put(-23,-22){\line(1,-1){19}}
\end{picture}
$$  
\end{remark}

\defin{}\label{def_type_admissible} 
Assume that $\{b_1,\ldots,b_s\}$ is some decomposition of a multilinear path $a$ in $\Q^{\ast}$ into primitive closed paths and a diagram $D$ is the type of $a$ with respect to $\{b_1,\ldots,b_s\}$. We say that $a$ is {\it admissible with respect to}  $\{b_1,\ldots,b_s\}$ if $D$ is admissible and the following conditions hold:
\begin{enumerate}
\item[a)]
If there is a triangle in $D$ with vertices $v_i,v_j,v_k$, then $b_i,b_j,b_k$ form a {\it fan}, i.e., there is a $u\in\Q_0$ such that $$\Ver{b_i}\cap\Ver{b_j}=\Ver{b_i}\cap\Ver{b_k}=\Ver{b_j}\cap\Ver{b_k}=\{u\}.$$

\item[b)]
If $v_i,v_j$ as well as $v_j,v_k$ are connected by means of edges marked with $2$, then $b_i,b_j,b_k$ form a {\it chain}, i.e., $b_j\sim c_1 c_2$ for paths $c_1,c_2$ with $c_1',c_1''\in\Ver{b_i}$ and $\Ver{c_1}\cap\Ver{b_k}=\emptyset$. Schematically, we depict this condition as follows:
$$ 
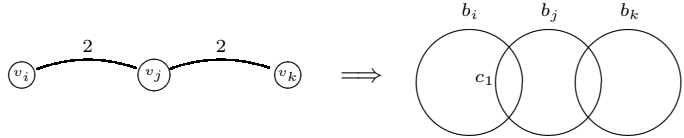
\xymatrix@C=1.3cm@R=1.3cm{ %
\vtx{v_i} \ar@/^/@{-}[r]^{2} & \vtx{v_j} \ar@/^/@{-}[r]^{2} & \vtx{v_k}  \\
}   \quad \Longrightarrow \qquad\qquad\qquad\qquad\qquad\qquad
\begin{picture}(0,30)(60,25)
\put(-30,25){\circle{50}} 
\put(0,25){\circle{50}} 
\put(30,25){\circle{50}} 
\put(-33,50){$\scriptstyle b_i$} 
\put(-3,50){$\scriptstyle b_j$} 
\put(27,50){$\scriptstyle b_k$} 
\put(-28,25){$\scriptstyle c_1$} 
\end{picture}
$$
\medskip

\noindent{}The orientations of closed paths $b_i,b_j,b_k$ can be arbitrary, so we do not specify it on the picture. In other words, for $\Ver{b_j}\cap\Ver{b_k}=\{u,v\}$ we do not have the following situation: 
$$ 
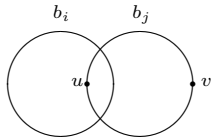
\begin{picture}(0,30)(0,25)
\put(-30,25){\circle{50}} 
\put(0,25){\circle{50}} 
\put(-20,25){\circle*{2}}\put(-26,24){$\scriptstyle u$} 
\put(20,25){\circle*{2}}\put(23,24){$\scriptstyle v$} 
\put(-33,50){$\scriptstyle b_i$} 
\put(-3,50){$\scriptstyle b_j$} 
\end{picture}
$$
\end{enumerate}
\bigskip

\noindent{}If $a$ is admissible with respect to some decomposition of $a$, then $a$ is called {\it admissible}; otherwise, we say that $a$ is not {\it admissible}.

\remark{} Obviously, using notations from the previous definition we have the following statement. Let paths $b_i,b_j$ intersect at a vertex $u$ ($i\neq j$) and $k\neq i,j$. Then $u\in\Ver{b_k}$ if and only if $b_i,b_j,b_k$ form a fan.   
\bigskip 

\remark{}\label{remark_no3}
If a multilinear path $a$ is admissible, then $\deg_v(a)\leq 3$ for all $v\in\Q_0$, where $\deg_v(a)$ is a number of closed primitive paths $a_1,\ldots,a_s$ with $a_1'=\cdots=a_s'=v$ such that $a\sim a_1\cdots a_s$.

\example\label{ex0} The path $h$ from Example~\ref{ex_decomp} is admissible with respect to decomposition from part~a) as well as part~b) of Example~\ref{ex_decomp} (see Example~\ref{ex_type}). 

\example\label{ex1} Let $\Q$ be the following quiver: 
$$\xymatrix@C=1.3cm@R=1.3cm{ %
\vtx{} \ar@/^/@{->}[r]^{x_3} \ar@/_/@{<-}[d]_{x_4} & %
\vtx{} \ar@/_/@{<-}[d]_{a_3} \ar@/^/@{->}[d] \ar@/_/@{->}[dr] & %
\vtx{} \ar@/_/@{<-}[d]_{b_2} \ar@/^/@{->}[d]^{b_1} & %
\vtx{} \ar@/^/@{<-}[d]^{c_3} \ar@/^/@{->}[dl]^{c_2} &\\
\vtx{} \ar@/_/@{<-}[r]_{x_1} & \vtx{} \ar@/_/@{<-}[r]_{a_1} & %
\vtx{} \ar@/_/@{->}[r]_{c_1} & \vtx{} \ar@/_/@{->}[r]_{y_1} \ar@/^/@{<-}[r]^{y_2} & \vtx{}\\
}
\begin{picture}(0,0)
\put(-120,-30){$\scriptstyle a_2$}%
\put(-135,-20){$\scriptstyle x_2$}%
\end{picture}
$$%
\noindent{}Consider the closed path $h=a_2 a_3 x_2 x_3 x_4 x_1 a_1 b_1 b_2 c_2 c_3 y_2 y_1 c_1$. Then $\{a,b,c,x,y\}$ is some decomposition of $h$ into primitive closed paths, where 
$a=a_1a_2a_3$, $b=b_1b_2$, $c=c_1c_2c_3$, $x=x_1x_2x_3x_4$, and $y=y_1y_2$. The type of $h$ with respect to $\{a,b,c,x,y\}$ is the following diagram $D$:
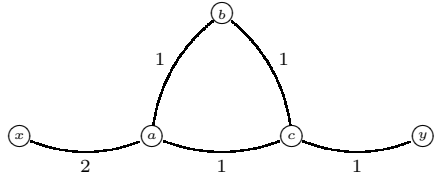
$$\xymatrix@C=0.6cm@R=1.3cm{ %
& & & \vtx{b} \ar@/_/@{-}[dl]_{1} \ar@/^/@{-}[dr]^{1}   & & &   \\ %
\vtx{x} \ar@/_/@{-}[rr]_{2} && \vtx{a} \ar@/_/@{-}[rr]_{1} && \vtx{c} \ar@/_/@{-}[rr]_{1} && \vtx{y} \\ %
}$$
\medskip%
\noindent{}Here vertices of $D$ are denoted by the corresponding closed primitive paths. By definition, $h$ is admissible with respect to $\{a,b,c,x,y\}$.

\example\label{ex1b} Let $\Q$ be the following quiver: 
$$
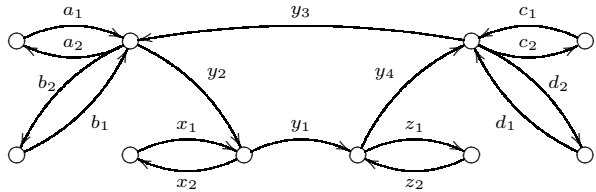
\xymatrix@C=1.3cm@R=1.3cm{ %
\vtx{} \ar@/^/@{->}[r]^{a_1} \ar@/_/@{<-}[r]^{a_2} & 
\vtx{} \ar@/^/@{<-}[rrr]^{y_3} \ar@/^/@{->}[rd]^{y_2} \ar@/_/@{->}[dl] \ar@/^/@{<-}[dl] &&&
\vtx{} \ar@/^/@{<-}[r]^{c_1} \ar@/_/@{->}[r]^{c_2} \ar@/^/@{->}[rd] \ar@/_/@{<-}[rd] &
\vtx{} \\
\vtx{} &
\vtx{} \ar@/^/@{->}[r]^{x_1} \ar@/_/@{<-}[r]_{x_2} &
\vtx{} \ar@/^/@{->}[r]^{y_1} &
\vtx{} \ar@/^/@{->}[ur]^{y_4} \ar@/^/@{->}[r]^{z_1} \ar@/_/@{<-}[r]_{z_2}&
\vtx{}& \vtx{}\\      
}
\begin{picture}(0,0)
\put(-190,-30){$\scriptstyle b_1$}%
\put(-210,-15){$\scriptstyle b_2$}%
\put(-37,-30){$\scriptstyle d_1$}%
\put(-17,-15){$\scriptstyle d_2$}%
\end{picture}$$%
\noindent{}For $a=a_1a_2$, $b=b_1b_2$, $c=c_1c_2$, $d=d_1d_2$, $x=x_1x_2$, $y=y_1y_2y_3y_4$, $z=z_1z_2$, we consider the closed path $h=y_1 x y_2 a b y_3 c d y_4 z$. Then $\{a,b,c,d,x,y,z\}$ is some decomposition of $h$ into primitive closed paths. The type of $h$ with respect to the given decomposition is the following diagram $D$:
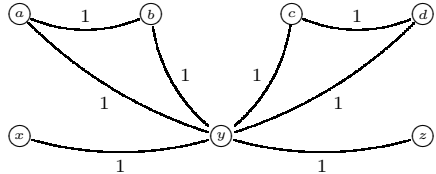
$$\xymatrix@C=0.6cm@R=1.3cm{ %
\vtx{a} \ar@/_/@{-}[rr]^{1} \ar@/_/@{-}[rrrd]_{1} && \vtx{b} \ar@/_/@{-}[rd]^{1} &&
\vtx{c} \ar@/_/@{-}[rr]^{1} \ar@/^/@{-}[ld]_{1}   && \vtx{d} \ar@/^/@{-}[llld]^{1}  \\ %
\vtx{x} \ar@/_/@{-}[rrr]_{1} &&& \vtx{y}\ar@/_/@{-}[rrr]_{1} &&&\vtx{z}\\
}$$%
\noindent{}Here vertices of $D$ are denoted by the corresponding closed primitive paths. By definition, $h$ is admissible with respect to $\{a,b,c,d,x,y,z\}$.

\defin{}\label{def_S1}  Denote by $\SetI$ a maximal (by inclusion) subset of admissible tree paths such that elements of $\SetI$ have pairwise different multidegrees.

\begin{theo}\label{theo_mgs_not2}
If $\Char{\FF}\neq2$, then $\{\det(X_a), \tr(X_b)\,|\, a\in\Q_1,\;b\in \SetI\}$ is a minimal generating set for $SI(\Q)$.
\end{theo}

\section{Generating set}\label{section_generating_set}

In this section we prove Theorem~\ref{theo_Fedotov} over a field of arbitrary characteristic. We have already shown that elements from Theorem~\ref{theo_Fedotov} belong to $SI(\Q)$.

The description of generators for semi-invariants of a quiver from~\cite{DZ01} was reformulated in Theorem~1 from~\cite{Lopatin_so_inv}, where more general notion of semi-invariants of supermixed representations of a quiver was considered. In the mentioned theorem semi-invariants were described using the notion of {\it tableau with substitution} $(T, (Y_1,\ldots,Y_l))$ and {\it block partial linearization of the pfaffian} ${\rm bpf}_T(Y_1,\ldots,Y_l)$ that were given in~\cite{Lopatin_bplp}.  In this article we only use a partial case of the notion of tableau with substitution.

\defin{} Assume that $n\in\NN$ and $Y_1,\ldots,Y_l$ are $n\times n$ matrices. Let $m=2l/n\in\NN$. A pair $(T,(Y_1,\ldots,Y_l))$ is called a {\it multilinear tableau with substitution} (m.t.s.) if
\begin{enumerate}
\item[$\bullet$] $T$ is an $n\times m$ tableau filled with arrows $\{1,\ldots,l\}$; 

\item[$\bullet$] {\it an arrow} goes from one cell of the tableau into another one, and each cell of the tableau is either the head or the tail of one and only one arrow.  
\end{enumerate}
We refer to $T$ as {\it tableau} of dimension $(n,\ldots,n)$ ($m$ times), and we write $\ga\in T$ for an arrow $\ga$ from $T$. Given an
arrow $\ga\in T$, denote by $\ga'$ and $\ga''$ the columns containing the head and the tail of $\ga$, respectively. Similarly, denote by $'\ga$ and $''\ga$ the rows containing 
the head and the tail of $\ga$, respectively. Schematically this is depicted as
$$
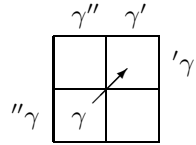
\begin{picture}(0,50)
\put(-10,30){%
\put(0,0){\rectangle{10}{10}\put(-3,-3){}}%
\put(20,0){\rectangle{10}{10}\put(-3,-3){}}%
\put(0,-20){\rectangle{10}{10}\put(-3,-3){$\ga$}\put(5,5){\vector(1,1){13}}}%
\put(20,-20){\rectangle{10}{10}\put(-3,-3){}}%
\put(-3,15){$\ga''$}%
\put(17,15){$\ga'$}%
\put(35,-3){$'\ga$}%
\put(-26,-23){$''\ga$}%
}%
\end{picture}
$$%

\example\label{ex_new} Let $T$ be the tableau
$$
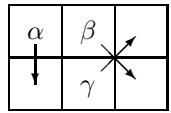
\begin{picture}(0,55)
\put(-20,40){%
\put(0,0){\rectangle{10}{10}\put(-3,-3){$\al$}\put(0,-5){\vector(0,-1){15}}}%
\put(20,0){\rectangle{10}{10}\put(-3,-3){$\be$}\put(5,-5){\vector(1,-1){13}}}%
\put(40,0){\rectangle{10}{10}}%
\put(0,-20){\rectangle{10}{10}}%
\put(20,-20){\rectangle{10}{10}\put(-3,-3){$\ga$}\put(5,5){\vector(1,1){13}}}%
\put(40,-20){\rectangle{10}{10}}%
}%
\end{picture}
$$
of dimension $(2,2,2)$ and $Y_1,Y_2,Y_3$ be $2\times 2$ matrices. Then $l=m=3$, $n=2$ and $(T,(Y_1,Y_2,Y_3))$ is an m.t.s. Note that $'\al=2$, $''\al=1$, and $\al'=\al''=1$. 

\bigskip
\defin{} Let $(T,(Y_1,\ldots,Y_l))$ be an m.t.s.~of dimension
$(n,\ldots,n)$ and $m=2s/n$. Define 
\begin{eq}\label{eq_def_F0}
{\rm bpf}_T(Y_1,\ldots,Y_l)=\sum_{\pi_1,\ldots,\pi_m\in S_{n}}
\sign(\pi_1)\cdots\sign(\pi_m)\prod_{\ga\in T}
(Y_\ga)_{\pi_{\ga''}(''\ga),\pi_{\ga'}('\ga)},
\end{eq}%
where $(Y_{\ga})_{ij}$ stands for the $(i,j)^{\rm th}$ entry of $Y_{\ga}$.
\bigskip
 
We assume that $\Q_0=\{1,\ldots,r\}$.    Theorem~1 of~\cite{Lopatin_so_inv} immediately implies that the algebra $SI(\Q)$ is generated by 
\begin{enumerate}
\item[a)] $\det(X_a)$, where $a\in\Q_1$;       

\item[b)] $\si_t(X_b)$, where $b$ is a closed path in $\Q$ and $t=1,2$;       

\item[c)] ${\rm bpf}_T(Y_1,\ldots,Y_l)$, where 
\begin{enumerate}
\item[$\bullet$] $(T,(Y_1,\ldots,Y_l))$ is an m.t.s.~of dimension $(2,\ldots,2)$ with arrows $\{\gamma_1,\ldots,\gamma_l\}$; 

\item[$\bullet$] $T$ is a union of $2r$ rectangular (possibly empty) blocks $A_1,\ldots,A_r$, $B_1,\ldots,B_r$ with two rows such that every cell of $T$ belongs to one and only one block;

\item[$\bullet$] every arrow $\gamma_k$ of the tableau $T$ goes from $A_i$ to $B_j$ for some $i,j$; moreover, there is a path $c_k$ in $\Q$ such that $c_k'=i$, $c_k''=j$, and $X_{c_k}=Y_k$.
 \end{enumerate}
\end{enumerate}
Note that $l$ is even. Obviously, a permutation of columns of $T$ does not affect ${\rm bpf}_T(Y_1,\ldots,Y_l)$ and a permutation of cells from a fixed column of $T$ changes ${\rm bpf}_T(Y_1,\ldots,Y_l)$ by $\pm1$. Hence without loss of generality we can assume that $T$ is equal to the following tableau $T_l$:
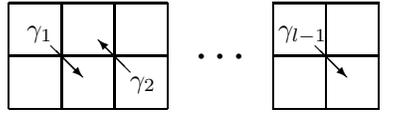
$$\begin{picture}(70,50)%
\put(-35,35){
\put(0,0){\rectangle{10}{10}\put(-3,-3){$\ga_1$}\put(6,-6){\vector(1,-1){12}}}%
\put(20,0){\rectangle{10}{10}\put(16,-22){$\ga_2$}\put(16,-16){\vector(-1,1){12}}}%
\put(40,0){\rectangle{10}{10}}%
\put(0,-20){\rectangle{10}{10}}%
\put(20,-20){\rectangle{10}{10}}%
\put(40,-20){\rectangle{10}{10}}%
\put(100,0){\rectangle{10}{10}\put(-8,-3){$\ga_{l-1}$}\put(6,-6){\vector(1,-1){12}}}%
\put(120,0){\rectangle{10}{10}}%
\put(100,-20){\rectangle{10}{10}}%
\put(120,-20){\rectangle{10}{10}}%
\put(70,-10){
\put(0,0){\circle*{2}} %
\put(-7,0){\circle*{2}} %
\put(7,0){\circle*{2}} %
}
\put(145,-10){,}
}
\end{picture}$$
where we have not depicted the arrow $\ga_l$ that goes from the bottom left cell to the top right cell.

\begin{lemma}\label{lemma_tabl}
Given $2\times 2$ matrices $Z_1,\ldots,Z_l$ $(l>0)$ over a commutative $\FF$-algebra, we have
$${\rm bpf}_{T_l}(Z_1,\ldots,Z_l)=(-1)^l\tr(Z_1 J_2 Z_2^T J_2 \cdots Z_{l-1} J_2 Z_l^T J_2).$$
\end{lemma}  
\begin{proof} We set $C_k=Z_k$ for odd $k$ and $C_k=Z_k^T$ for even $k$.
By definition, 
$${\rm bpf}_{T_l}(Z_1,\ldots,Z_l)=\sum_{\pi_1\in
S_2,\ldots,\pi_l\in S_2}
\sign(\pi_1)\cdots\sign(\pi_l)\prod_{k=1}^l
(C_k)_{\pi_{k}(1),\pi_{k+1}(2)},
$$
where we assume that $\pi_{l+1}=\pi_1$. On the other hand, 
$$\tr(C_1 J_2\cdots C_l J_2)=\sum_{1\leq i_1,\ldots,i_l\leq 2} \prod_{k=1}^l(C_k J_2)_{i_k,i_{k+1}},$$%
where $i_{l+1}=i_1$. Let $\tau_k\in S_2$ satisfies $\tau_k(1)=i_k$. Then $\tau_k(2)=\xi(i_{k})$ and $\sign(\tau_k)=-(-1)^{i_k}$, where $\xi$ is the non-identical permutation from $S_2$. The fact that $(C_k J_2)_{ij}=(-1)^j (C_k)_{i,\xi(j)}$ completes the proof. 
\end{proof}

Since $c_k''=c_{k+1}''$ for odd $k$ and $c_k'=c_{k+1}'$ for even $k$ ($1\leq k\leq l$), where $c_{l+1}$ stands for $c_1$, $e=c_1c_2^{\ast}\ldots c_{l-1} c_l^{\ast}$ is a closed path in $\Q^{\ast}$. Lemma~\ref{lemma_tabl} implies that ${\rm bpf}_T(Y_1,\ldots,Y_l)=\pm \tr(X_e)$. Relation~(D) from Lemma~\ref{lemma_relations} (see below) completes the proof of Theorem~\ref{theo_Fedotov}.

\section{Some relations}\label{section_relations}

In what follows, we write $\si_t(a)$ for $\si_t(X_a)$, where $a$ is a closed path $\Q^{\ast}$. 

For a $v\in\Q_0$ we denote by $1_v$ the {\it empty path} in the vertex $v$. We set $\Ver{1_v}=\{v\}$ and $\Arr{1_v}=\emptyset$. Given a path
$a$ with $a'=v$, we assume $1_v a=a$ and for a path $a$ with $a''=v$ we assume $a 1_v=a$. 
Denote by $\path{\Q^{\ast}}$ the set of all paths and empty paths in $\Q^{\ast}$. 

\begin{lemma}\label{lemma_relations} For closed paths $a,b,c$ and paths $x,x_1,x_2,y_1,y_2$ in $\Q^{\ast}$ the following relations hold.
\begin{enumerate}
 \item[(0)] $\si_t(a^{\ast})=\si_t(a)$,  $\si_t(y_1 y_2)=\si_t(y_2 y_1)$, where $t=1,2$; $\det(a b)\equiv0$.

 \item[(A)] $\tr(a^2 b)\equiv0$, where $a'=b'$, i.e.,
\begin{picture}(0,0)(-30,5)%
$\loopR{0}{0}{a\,}
\xymatrix@C=1.3cm@R=1.3cm{\vtx{}}
\loopL{0}{0}{\,b}\quad.$%
\end{picture}%
\medskip
 
 \item[($A'\!$)] $\tr(a b c)\equiv -\tr(acb)$, where $a'=b'=c'$.

 \item[(B)] $\tr(a^{\ast}b)\equiv -\tr(a b)$, where $a'=b'$.

 \item[(C)] $\tr(x x^{\ast} a)\equiv0 $, where $x'=a'$, i.e.,
\begin{picture}(0,0)(-20,5)
$\loopR{0}{0}{a} %
\xymatrix@C=1.3cm@R=1.3cm{ \vtx{}%
\ar@/^/@{<-}[r]^{x} & \vtx{}\\}%
\quad.$
\end{picture}
\medskip

 \item[($C'\!$)] $\tr(x_1 x_2^{\ast} a)\equiv -\tr(x_2 x_1^{\ast} a)$, where $x_1'=x_2'=a'$.

 \item[(D)] $\tr(x^{\ast} y_1 x^{\ast} y_2)\equiv0$, where $x'=y_1'=y_2'$ and $x''=y_1''=y_2''$, i.e.,
\begin{picture}(0,0)(-20,5)
$\xymatrix@C=1.3cm@R=1.3cm{ %
\vtx{}\ar@/^/@{<-}[r]^{x^{\ast}} \ar@2@/_/@{->}[r]_{y_1,y_2} & \vtx{}\\}%
\quad.$
\end{picture}
\bigskip

 \item[($D'\!$)] $\tr(x_1^{\ast} y_1 x_2^{\ast} y_2)\equiv -\tr(x_2^{\ast} y_1 x_1^{\ast} y_2)$, where $x_1'=x_2'=y_1'=y_2'$ and $x_1''=x_2''=y_1''=y_2''$.

 \item[(E)] $\tr(x x^{\ast})=2\det(x)$.
\end{enumerate}
\end{lemma}
\begin{proof}
Relations (0) and (E) are trivial. Relation~(A) follows from 
$$\tr(a^2 b)=\tr(a) \tr(a b) - \det(a)\tr(b).$$

Relation~(B) follows from
$$\tr(ab)=-\tr(a^{\ast}b) + \tr(a)\tr(b).$$

Relation~(C) follows from
$$\tr(x x^{\ast} a)=\det(x)\tr(a).$$

Relation~(D) follows from
$$\tr(x^{\ast} y_1 x^{\ast} y_2) = - \det(x) \tr(y_1^{\ast} y_2) + \tr(x^{\ast} y_1)\tr(x^{\ast} y_2).$$

Applying linearization to (C) and (E), i.e., making a substitution $X_x\to X_{x_1} + X_{x_2}$, where $x_1'=x_2'=x'$ and $x_1''=x_2''=x''$, and taking the homogeneous component of degree $1$ with respect to both $X_{x_1}$ and $X_{x_2}$, we obtain relations ($\rm C'$) and ($\rm D'$). The proof of ($\rm A'$) is similar. 
\end{proof}

\begin{lemma}\label{lemma_two_rel}
Let $\Char{\FF}\neq2$ and $a,b,c,d$ be paths in $\Q^{\ast}$ that are closed in $v\in \Q_0$. Then
\begin{enumerate}
\item[($R_1\!$)] $\tr(abcd)\equiv0$;

\item[($R_2\!$)] $\tr(abc)\equiv0$, where $a$ and $b$ intersect at a vertex different from $v$.
\end{enumerate}
\end{lemma}
\begin{proof}
Applying ($\rm A'$) several times, we obtain the proof of ($\rm R_1$):
$$\tr(a\cdot b\cdot cd)\equiv - \tr(ac\cdot d\cdot b)\equiv \tr(a c b d)=\tr(da\cdot c\cdot b)\equiv -\tr(dabc).$$
Assume that $a=x_1 y_1$ and $b=x_2 y_2$ for paths $x_1,x_2,y_1,y_2$ in $\Q^{\ast}$ with $x_1''=x_2''=y_1'=y_2'=w$ for a vertex $w$ and $w\neq v$. By ($\rm A'$) we have 
$$\tr(abc)=\tr(c\cdot x_1y_1 \cdot x_2y_2)\equiv -\tr(c x_2y_2 x_1y_1).$$
Applying ($\rm D'$) two times we obtain
$$\tr(c x_2y_2 x_1y_1)\equiv -\tr(c x_2y_1 x_1y_2)\equiv\tr(c x_1y_1 x_2y_2)=\tr(cab). $$ 
Thus, ($\rm R_2$) is proven.
\end{proof}

Note that if $a\sim b$ for closed paths $a$ and $b$, then $\mdeg(a)=\mdeg(b)$ and $\tr(a)=\tr(b)$. The next lemma generalizes this remark. 

\begin{lemma}\label{lemma_mdeg_tr}
Assume that $a,b$ are closed paths in $\Q^{\ast}$ and $\mdeg(a)=\mdeg(b)$. Then $\tr(a)\equiv \pm\tr(b)$.
\end{lemma}
\begin{proof} We assume $a=a_1\cdots a_s$, where $a_i\in\Q_1^{\ast}$. Since $\mdeg(a)=\mdeg(b)$, the equivalence $b\sim a_1 b_2\cdots b_s$ holds, where $b_i\in\Q_1^{\ast}$. If $s=1$, then $b\sim a$.

Let $s>1$. Since $\mdeg(a_2\cdots a_s)=\mdeg(b_2\cdots b_s)$, we have $b_2\cdots b_s=c a_2 d$ or $b_2\cdots b_s=c a_2^{\ast} d$ for $c,d\in\path{\Q^{\ast}}$. By case by case consideration we will show that $\tr(b)\equiv\pm\tr(a_1a_2e)$ for $e\in\path{\Q^{\ast}}$. Repeating this procedure we complete the proof.

Case 1). Assume that $c$ is not empty and $a_2$ is not a loop.

Case 1.1). Let $b_2 \cdots b_s = c a_2 d$. Then we depict the closed path $a_1 b_2\cdots b_s$ in $\Q^{\ast}$ as follows:
$$
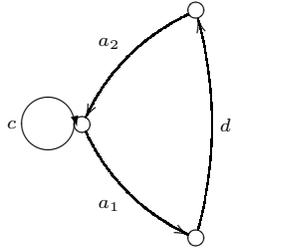
\loopR{0}{43}{c}%
\xymatrix@C=1.3cm@R=1.3cm{ %
&  \vtx{}    \\
\vtx{} \ar@/^/@{<-}[ru]^{a_2} \ar@/_/@{->}[rd]_{a_1}&    \\
&  \vtx{}\ar@/_/@{->}[uu]_{d}    \\
}%
\begin{picture}(0,0)(0,0)
\put(20,-45){,} 
\end{picture}
$$
where $d$ can be empty. Since $\mdeg(a)=\mdeg(a_1 a_2 d) + \mdeg(c)$ and $a=a_1a_2\cdots a_s$, we have that $\Ver{c}\cap \Ver{d}$ is not empty. Denote $v=a_1''$.

Case 1.1.a). Let $v\in\Ver{c}\cap \Ver{d}$.  Then $d=d_1 d_2$ for a path $d_1$ in $\Q^{\ast}$ with $d_1''=v$ and $d_2\in\path{\Q^{\ast}}$. Hence we have
$$
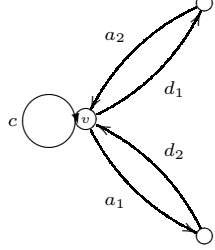
\loopR{-1}{45}{c}%
\xymatrix@C=1.3cm@R=1.3cm{ %
&  \vtx{}    \\
\vtx{v} \ar@/^/@{<-}[ru]^{a_2} \ar@/_/@{->}[ru]_{d_1} 
        \ar@/_/@{->}[rd]_{a_1} \ar@/^/@{<-}[rd]^{d_2}&    \\
&  \vtx{}    \\
}%
\begin{picture}(0,0)(0,0)
\put(20,-45){.} 
\end{picture}
$$
Applying relation~($\rm A'$), we obtain $\tr(b)=\tr(d_2a_1\cdot c\cdot a_2d_1)\equiv -\tr(a_1a_2d_1 c d_2)$.

Case 1.1.b). Let $w\in\Ver{c}\cap \Ver{d}$ for a vertex $w$ with $v\neq w$.  Then $c=c_1 c_2$ and $d=d_1 d_2$  for $d_1,d_2\in\path{\Q^{\ast}}$ and paths $c_1,c_2$ with $c_2'=d_2'=w$. Hence we have
$$
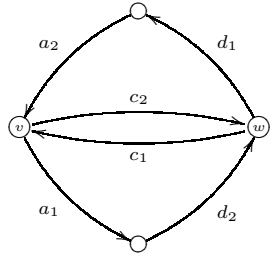
\xymatrix@C=1.3cm@R=1.3cm{ %
&  \vtx{}   &   \\
\vtx{v} \ar@/^/@{<-}[ru]^{a_2} \ar@/_/@{->}[rd]_{a_1} 
        \ar@/^/@{->}[rr]^{c_2} \ar@/_/@{<-}[rr]_{c_1}&&   
\vtx{w}  \ar@/_/@{->}[ul]_{d_1} \ar@/^/@{<-}[dl]^{d_2} 
 \\
&  \vtx{} &    \\
}%
\begin{picture}(0,0)(0,0)
\put(20,-45){.} 
\end{picture}
$$
Applying relation~($\rm D'$), we obtain $\tr(b)=\tr(c_1\cdot c_2\cdot a_2 d_1\cdot d_2 a_1)\equiv -\tr(a_1 a_2 d_1 c_2 c_1 d_2)$.

Case 1.2). Let $b_2\cdots b_s=c a_2^{\ast} d$. Thus we depict $a_1 b_2\cdots b_s$ as follows:
$$
\xymatrix@C=1.3cm@R=1.3cm{ %
&  \vtx{}    \\
\vtx{} \ar@/^/@{->}[ru]^{a_2^{\ast}} \ar@/_/@{<-}[ru]_{c} 
        \ar@/_/@{->}[rd]_{a_1} \ar@/^/@{<-}[rd]^{d}    \\
&  \vtx{}    \\
}%
\begin{picture}(0,0)(0,0)
\put(20,-45){.} 
\end{picture}
$$
By relation (B), $\tr(b)=\tr(c a_2^{\ast}\cdot d a_1)\equiv -\tr(a_2 c^{\ast} d a_1)=-\tr(a_1 a_2 c^{\ast} d)$.

Case 2). Let $c$ be a non-empty path and $a_2$ be a loop. 

If $b_2\cdots b_s = c a_2 d$, then relation~($\rm A'$) implies that $\tr(b)=\tr(d a_1 c a_2)\equiv - \tr(a_1 a_2 c d)$.

If $b_2\cdots b_s = c a_2^{\ast} d$, then relations~($\rm A'$) and (B) imply that $\tr(b)=\tr(d a_1 c a_2^{\ast})\equiv - \tr(a_1 a_2^{\ast} c d) \equiv \tr(a_1 a_2 c d)$.

Case 3). Let $c$ be empty.

If $b_2\cdots b_s = a_2 d$, then  $\tr(b)= \tr(a_1 a_2 d)$.

If $b_2\cdots b_s = a_2^{\ast} d$, then applying relation  (B) we obtain $\tr(b)= \tr(a_1 a_2^{\ast} d) \equiv -\tr(a_1 a_2 d)$.
  
Since we have considered all cases, the proof is completed.  
\end{proof}

\begin{lemma}\label{lemma_not_tree_path} 
Assume that $a$ is a path in $\Q^{\ast}$ such that $a$ is not a tree path and $a\not\sim x x^{\ast}$ for any $x\in\Q_1$. Then $\tr(a)\equiv0$.
\end{lemma}
\begin{proof} Since $a$ is not a tree path, then one of the following two cases holds.

Case 1). Let $\deg_x(a)\geq2$ for an $x\in\Q_1^{\ast}$. Then relations~(A), ($\rm A'$) and~(D) imply that $\tr(a)\equiv0$.

Case 2). Assume that there exists an $x\in\Q_1^{\ast}$ such that $a\sim a_1 x a_2 x^{\ast}$ for $a_1,a_2\in\path{\Q^{\ast}}$ satisfying one of the following conditions: 
\begin{enumerate}
\item[$\bullet$] $a_1$ or $a_2$ is empty; 

\item[$\bullet$] there is a $w\in\Q_0$ such that $w\in\Ver{a_1}\cap\Ver{a_2}$. 
\end{enumerate}
If $a_1$ or $a_2$ is empty, then $\tr(a)\equiv0$ by relation~(C).  

Assume that $a_1$ and $a_2$ are not empty. If $x$ is a loop, then relations ($\rm A'$) and (C) imply the required equality. 

Assume that $x$ is not a loop. Then $a_i=c_i d_i$ for $c_i,d_i\in\path{\Q^{\ast}}$, $i=1,2$,  satisfying $c_1''=c_2''=w$. Denote $x'=u$ and $x''=v$. We can depict $a_1 x a_2 x^{\ast}$ as follows:
$$ 
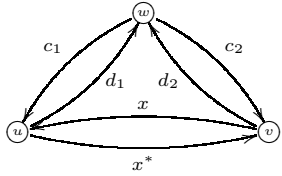
\xymatrix@C=1.3cm@R=1.3cm{ %
&\vtx{w} %
\ar@/_/@{->}[ld]_{c_1} \ar@/^/@{<-}[ld]
\ar@/^/@{->}[rd]^{c_2} \ar@/_/@{<-}[rd] & \\
\vtx{u} \ar@/_/@{->}[rr]_{x^{\ast}} \ar@/^/@{<-}[rr]^{x} &&\vtx{v}  }%
\begin{picture}(0,0)(0,0)
\put(-67,-25){$\scriptstyle d_1$} 
\put(-47,-25){$\scriptstyle d_2$} 
\put(10,-25){.} 
\end{picture}
$$

Case 2.1). If there is an $i=1,2$ such that $c_i$ or $d_i$ is empty, then $w\in\{u,v\}$. If $w=u$, then relations ($\rm A'$) and (C) imply that %
$\tr(a)=\tr(a_1\cdot x c_2 \cdot d_2 x^{\ast})\equiv -\tr(a_1 d_2 x^{\ast} x c_2)\equiv 0$. If $w=v$, then we obtain the required equality similarly.    

Case 2.2). If $c_1,c_2,d_1,d_2$ are non-empty paths, then applying ($\rm D'$) and (C) we obtain %
$\tr(a)=\tr(c_1 d_1 x c_2 d_2 x^{\ast})= 
\tr(d_1 x \cdot c_2 \cdot d_2 \cdot x^{\ast} c_1)\equiv 
-\tr(d_1 x \cdot x^{\ast} c_1 \cdot d_2 \cdot c_2) \equiv0$. The proof is completed. 
\end{proof}

\section{The case of $\Char{F}=2$}\label{section_case2}

In this section we assume that $\Char{\FF}=2$ unless otherwise stated. We say that a tree path $a$ in $\Q^{\ast}$ is {\it simple} if for every $x\in\Arr{a}$ with  $x^{\ast}\not\in\Arr{a}$ we have that $x$ is a loop. In other words, every block of simple tree path is a quiver with one vertex and several loops. We use the following remark in the next two sections.

\remark{}\label{remark_Phi}
To define a homomorphism $\Phi: SI(\Q)\to R$ of $\FF$-algebras, where $R$ is a commutative $\FF$-algebra, for every $z\in\Q_1$ we will specify $2\times 2$ matrix $\Phi(X_z)$ over $R$. Then we set that $\Phi(x_{ij}^{z})$ is the $(i,j)^{\rm th}$ entry of $\Phi(X_z)$. Note that in some cases we define only $\Phi(X_{z^{\ast}})$, not $\Phi(X_{z})$. Then we assume that 
\begin{eq}\label{eq_astast}
\Phi(X_z)=\Phi(X_{z^{\ast}})^{\ast}.
\end{eq}%
With abuse of notation in some cases we define $\Phi(X_z)$ together with $\Phi(X_{z^{\ast}})$. In these cases the equality~\Ref{eq_astast} holds.

\begin{lemma}\label{lemma_simple_tree_path}
If $a$ is a simple tree path in $\Q^{\ast}$, then $\tr(a)$ is indecomposable.
\end{lemma}  
\begin{proof}
Let $\Ver{a}=\{v_1,\ldots,v_l\}$. By definition of simple tree path, we have  $$\Arr{a}=\{a_1,a_1^{\ast},\ldots,a_{l-1},a_{l-1}^{\ast}\}\bigcup \bigcup_{i=1}^l\{b_{i,1},\ldots,b_{i,t_i}\},$$
where $t_1,\ldots,t_l\geq0$ and $b_{i,1},\ldots,b_{i,t_i}$ are loops in $v_i$ ($1\leq i\leq l$) such that $b_{i,j}\not\sim b_{i,k}$ for $j\neq k$. 

We prove the lemma by induction on $l\geq 1$. 

Let $l=1$.  If $t_1=1$, then $\tr(a)=\tr(b_{11})\not\equiv0$. 

Assume that $t_1\geq2$. Define a homomorphism $\Phi: SI(\Q)\to \FF[x_{ij}^{b_{11}}\,|\,1\leq i,j\leq 2]$ of $\FF$-algebras as follows: for every $y\in\Q_1^{\ast}$ we set
$$\Phi(X_y)=\left\{
\begin{array}{rl}
X_y,& \text{if } y=b_{11}\\
E,& \text{otherwise }\\
\end{array}
\right.
$$
(see Remark~\ref{remark_Phi}). If $\tr(a)\equiv 0$, then $\tr(a)=\sum_q \al_q f_q h_q$ for $\al_q\in\FF$ and some products of traces $f_q,h_q$.  We have $\Phi(\tr(a))=x_{11}^{b_{11}}+x_{22}^{b_{11}}\neq0$. On the other hand, the equality  $\tr(E)=0$ implies $\Phi(\tr(a))=0$; a contradiction.

Induction step. Let $l\geq2$ and $\tr(a)$ be decomposable, i.e., $\tr(a)=\sum_q \al_q f_q$ for $\al_q\in\FF$ and some products of traces and determinants $f_q$ with two or more factors. Applying relation (B), we assume that there is no  $\tr(b_{ij}^{\ast}c)$, where $c\in\path{\Q^{\ast}}$, among traces in $f_q$.   Without loss of generality we can assume that $v_1$ corresponds to a leaf of the tree of $a$ (see the definition of tree path for details). Hence $t_1>0$. Moreover, without loss of generality we can assume that $a_1$ connects $v_1$ and $v_2$. Therefore,  $v_1\not\in\{a_i',a_i''\}$ for all $2\leq i\leq l-1$. 

Let $k=2$ satisfy the following property:
\begin{eq}\label{eq_tt}
t_k>0 \text{ or } v_k\in\{a_i',a_i''\} \text{ for at least three pairwise different $i$ with } 1\leq i\leq  l-1. 
\end{eq}%
Denote 
$$\Omega=\{a_2,a_2^{\ast},\ldots,a_{l-1},a_{l-1}^{\ast}\}\bigcup \bigcup_{i=2}^{l}\{b_{i,1},\ldots,b_{i,t_i}\}.$$
Then $\Arr{a}=\{a_{1},a_{1}^{\ast}\}\cup\{b_{1,1},\ldots,b_{1,t_1}\}\cup\Omega$. Schematically this is depicted as
$$\loopR{0}{0}{b_{11},\ldots,b_{1t_1}}
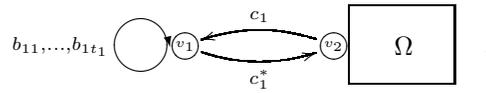
\xymatrix@C=1.5cm@R=1.5cm{ %
\vtx{v_1}\ar@/^/@{<-}[r]^{c_{1}}  \ar@/_/@{->}[r]_{c_{1}^{\ast}} &\vtx{v_2}} %
\begin{picture}(0,0)
\put(19,3){\rectangle{20}{15}\put(-3,-4){$\Omega$}}%
\put(50,0){,}
\end{picture}
$$%
where $c_1$ stands for $a_1$ or $a_{1}^{\ast}$. Define a homomorphism $\Psi: SI(\Q)\to \FF[x_{ij}^{y}\,|\,1\leq i,j\leq 2,\; y\in\Q_1]$ of $\FF$-algebras as follows: for every $y\in\Q_1^{\ast}$ we set 
$$\Psi(X_y)=\left\{
\begin{array}{rl}
E,& \text{if }  y\in\Arr{a}\backslash\Omega\\
X_y,& \text{otherwise}\\
\end{array}
\right.
$$
Let $c$ be a closed path in $\Q^{\ast}$ with $\Arr{c}\subset \Arr{a}$. Then
\begin{enumerate}
\item[$\bullet$] if $\deg_x(c)\geq1$ for an $x\in\Omega$, then $\Psi(\tr(c))=\tr(d)$ for a closed path $d$ with $\Arr{d}\subset \Omega$;

\item[$\bullet$] if $\Arr{c}\cap\,\Omega$ is empty, then $\Psi(\tr(c))=\tr(E)=0$;

\item[$\bullet$] $\Psi(\det(a_i))=\det(a_i)$ for all $2\leq i\leq l-1$;

\item[$\bullet$] if $\det(a_1)$ is a factor of $f_q$, then $\tr(b_{1,j_1}\cdots b_{1,j_s})$ is also a factor of $f_q$ for some $j_1,\ldots,j_s$; thus,  $\Psi(f_q)=0$.
\end{enumerate}
This remark implies that $\Psi(\tr(a))=\tr(e)\equiv0$ for a  path $e$ with $\Arr{e}= \Omega$. By condition~\Ref{eq_tt}, $e$ is a simple tree path. Since the tree of $e$ has exactly $l-1$ vertices, the induction hypothesis implies a contradiction. 

Let $k=2$ do not satisfy property~\Ref{eq_tt}. Then without loss of generality we can assume that there is a $k>1$ satisfying property~\Ref{eq_tt} such that $2,\ldots,k-1$ do not satisfy property~\Ref{eq_tt} and $\{a_{i-1}',a_{i-1}''\}\cap\{a_i',a_i''\}\neq\emptyset$ for all $2\leq i\leq k-1$. Schematically this is depicted as
$$\loopR{0}{0}{b_{11},\ldots,b_{1t_1}}
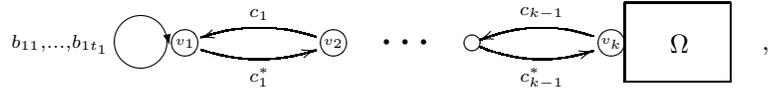
\xymatrix@C=1.5cm@R=1.5cm{ %
\vtx{v_1}\ar@/^/@{<-}[r]^{c_1}  \ar@/_/@{->}[r]_{c_1^{\ast}} &\vtx{v_2} %
& \vtx{} %
\ar@/^/@{<-}[r]^{c_{k-1}}  \ar@/_/@{->}[r]_{c_{k-1}^{\ast}} &\vtx{v_k}} %
\begin{picture}(0,0)
\put(-85,3){
\put(0,0){\circle*{2}} %
\put(-7,0){\circle*{2}} %
\put(7,0){\circle*{2}} %
}
\put(18,3){\rectangle{20}{15}\put(-3,-4){$\Omega$}}%
\put(50,0){,}
\end{picture}
$$%
where $c_i$ stands for $a_i$ or $a_i^{\ast}$ ($1\leq i\leq k-1$) and  
$$\Omega=\{a_k,a_k^{\ast},\ldots,a_{l-1},a_{l-1}^{\ast}\}\bigcup \bigcup_{i=k}^{l}\{b_{i,1},\ldots,b_{i,t_i}\}.$$%
Repeating the above reasoning we obtain a contradiction with the induction hypothesis. 
\end{proof}

\remark{}\label{remark_glue}
Let $u$ and $v$ be two different vertices of $\Q$. Denote by $\Q_{uv}$ the quiver that is the result of gluing of $u$ with $v$, i.e., $(\Q_{uv})_0=\Q_0\backslash\{v\}$ and $(\Q_{uv})_1=\{\tilde{x}\,|\,x\in\Q_1\}$, where
$$\tilde{x}'=\left\{
\begin{array}{rl}
x',& \text{if }  x'\neq v\\
u,& \text{otherwise}\\
\end{array}
\right.
\quad\text{ and }\quad
\tilde{x}''=\left\{
\begin{array}{rl}
x'',& \text{if }  x''\neq v\\
u,& \text{otherwise}\\
\end{array}
\right..
$$
Then $(\Q^{\ast})_{uv} = \Q^{\ast}_{uv}$. Let $a$ be a closed path in $\Q^{\ast}$ and let $b$ be the image of $a$ in $\Q^{\ast}_{uv}$. If $\tr(a)\equiv0$ in $SI(\Q)$, then $\tr(b)\equiv0$ in $SI(\Q_{uv})$. Similar result is valid for gluing of several vertices. 
\bigskip

The following lemma generalizes Lemma~\ref{lemma_simple_tree_path} for the case of arbitrary tree path.

\begin{lemma}\label{lemma_tree_path}
If $a$ is a tree path in $\Q^{\ast}$, then $\tr(a)$ is indecomposable.
\end{lemma}
\begin{proof} Assume that $\tr(a)\equiv0$. Let $\Gamma^{(1)},\ldots,\Gamma^{(l)}$ be blocks of the tree path $a$ (see Definition~\ref{def_tree_path}). 

For every $i$ we glue all vertices of $\Gamma^{(i)}$ together and denote the resulting quiver by $\Lambda$ (see Remark~\ref{remark_glue} for details). Let $b$ be the image of $a$ in $\Lambda$. Then Remark~\ref{remark_glue} implies that $\tr(b)\equiv0$ in $SI(\Lambda)$. On the other hand, $b$ is a simple tree path in $\Lambda$; a contradiction to Lemma~\ref{lemma_simple_tree_path}. 
\end{proof}

\begin{proof_of}{of Theorem~\ref{theo_mgs_2}.}
Denote by $P$ the set from the formulation of the theorem. 
Let the characteristic of $\FF$ be arbitrary. Theorem~\ref{theo_Fedotov} together with Lemmas~\ref{lemma_mdeg_tr},~\ref{lemma_not_tree_path} and relation~(E) show that $P$ generates $SI(\Q)$. 

Let $\Char{\FF}=2$. Then Lemma~\ref{lemma_tree_path}, the indecomposability of $\det(a)$ for all $a\in\Q_1$ together with the fact that all elements of $P$ have pairwise different multidegrees imply that $P$ is a minimal generating set for $SI(\Q)$.
\end{proof_of}

\section{The case of $\Char{F}\neq2$}\label{section_case_not2}

In this section we assume that $\Char{F}\neq2$.

\begin{lemma}\label{lemma_multilinear_path}
Let $a$ be a closed path in $\Q^{\ast}$ and $b$ be a multilinear path in $\supp(a)$ satisfying the following condition:  
\begin{eq}\label{eq_cond_old_multilin}
\text{ for all } x\in\Q_1 \text{ we have }\deg_x(b)=0 \text{ or } \deg_{x^{\ast}}(b)=0.
\end{eq}
If $\tr(b)\equiv0$, then $\tr(a)\equiv0$.
\end{lemma}
\begin{proof} Let $b=b_1\cdots b_r$ for $b_i\in\Q_1^{\ast}$ and $\De=\mdeg(a)-\mdeg(b)\in\NN^{\#\Q_1}$. Consider $v=b_1''$. Let  $$P_{v,\De}=\{x\in\Q_1\,|\,v\in\Ver{x} \text{ and }\De_x>0\}$$ 
be a non-empty set. Then there is a path $c_1$ in $\Q^{\ast}$ such that $c_1'=c_1''=v$, $\De^{(1)}=\De-\mdeg(c_1)\in\NN^{\#\Q_1}$ and the set $P_{v,\De^{(1)}}$ is empty. Moreover, we assume that the degree of $c_1$ is maximal. If $P_{v,\De}$ is empty, then we set $c_1$ is the empty path in the vertex $v$ and $\De^{(1)}=\De$. Then apply this procedure to $b_2'',\De^{(1)}$ to obtain $c_2,\De^{(2)}$ and so on. Finally, we construct a closed path $c=b_1 c_1\cdots b_r c_r$ in $\Q^{\ast}$ with $\mdeg(a)=\mdeg(c)$, where $c_i$ is either an empty path or a closed path in $\Q^{\ast}$. By Lemma~\ref{lemma_mdeg_tr}, $\tr(a)\equiv\pm \tr(c)$. Since $\tr(b)\equiv0$, we have 
\begin{eq}\label{eq_1}
\tr(b)=\sum_q \al_q f_q h_q,
\end{eq}%
where $\al_q\in\FF$ and $f_q,h_q\in SI(\Q)$ are homogeneous of positive degree. We apply the substitution $b_i\to b_i c_i$ for all $1\leq i\leq r$ to~\Ref{eq_1}. Since $b_i\not\sim b_j$ for $i\neq j$, this substitution is well defined. As the result, we obtain that $\tr(c)$ is decomposable. Thus, $\tr(a)\equiv0$.
\end{proof}

\begin{lemma}\label{lemma_trAB0}
Let $c=a_1a_2a_3b_1b_2b_3$ be a closed path in $\Q^{\ast}$, where
$a_i,b_i$ are such paths that $a_i'=b_i'$ ($1\leq i\leq3$). Then $\tr(c)\equiv0$.
\end{lemma}
\begin{proof}
By relation ($\rm D'$), we have $\tr(c)=\tr(a_1 a_2\cdot a_3\cdot b_1 b_2 \cdot b_3)\equiv - \tr( b_1\cdot b_2 a_3\cdot a_1\cdot a_2 b_3) \equiv \tr(b_2\cdot a_3 b_1\cdot a_2\cdot b_3 a_1)\equiv -\tr(c)$.
\end{proof}

\begin{cor}\label{cor_trAB0}
Let $a$ be a multilinear path in $\Q^{\ast}$ such that $\tr(a)\not\equiv0$. Assume that some decomposition of $a$ into primitive closed paths contains $b$ and $c$ with $b\not\sim c$. Then $\#(\Ver{b}\cap\Ver{c})\leq 2$.
\end{cor}
\begin{proof} For every $x\in\Q_1$ with $\deg_x(a)=\deg_{x^{\ast}}(a)=1$ we add a new arrow $\tilde{x}$ to $\Q$ with $\tilde{x}'=x'$ and $\tilde{x}''=x''$ and substitute $\tilde{x}^{\ast}$ for $x^{\ast}$ in $a$. Let $\tilde{a}$ be the resulting multilinear path in the resulting quiver $\widetilde{\Q}^{\ast}$. Note that $\tr(\tilde{a})$ is indecomposable in $SI(\widetilde{\Q}^{\ast})$.  Therefore, without loss of generality we can assume that $a$ satisfies condition~\Ref{eq_cond_old_multilin}.

Let $\#(\Ver{b}\cap\Ver{c})\geq 3$. Then $b\sim b_1 b_2 b_3$ and $c\sim c_1 c_2 c_3$ for paths $b_i,c_i$ in $\Q^{\ast}$ with $b_i'=c_i'$ ($1\leq i\leq3$). Thus Lemma~\ref{lemma_trAB0} implies that  $\tr(e)\equiv0$ for $e=b_1 b_2 b_3 c_1 c_2 c_3$. Since $a$ is a multilinear path satisfying~\Ref{eq_cond_old_multilin}, then $e$ is also a multilinear path satisfying~\Ref{eq_cond_old_multilin}. By Lemma~\ref{lemma_multilinear_path}, we obtain a contradiction. 
\end{proof}

In the formulation of the next lemma we use notions from Definition~\ref{def_type_admissible}. 

\begin{lemma}\label{lemma_possibilities}
Let $a$ be a multilinear path in $\Q^{\ast}$ such that $\tr(a)\not\equiv0$. Assume that some decomposition of $a$ into primitive closed paths contains $b_1,b_2,b_3$ such that  $b_{i}\not\sim b_j$ for $i\neq j$. Then up to permutation of indices of $b_1,b_2,b_3$ one of the following possibilities holds:
\begin{enumerate}
\item[a)] $\Ver{b_i}\cap\Ver{b_3}$ is empty for $i=1,2$;

\item[b)] $\Ver{b_1}\cap\Ver{b_3}$ is empty and either $b_1,b_2,b_3$ form a chain or 
$$\#(\Ver{b_1}\cap\Ver{b_2})\leq2 \text{ and } \#(\Ver{b_2}\cap\Ver{b_3})=1;$$

\item[c)] paths $b_1,b_2,b_3$ form a fan.
\end{enumerate} 
\end{lemma}
\begin{proof} As in the proof of Corollary~\ref{cor_trAB0},  without loss of generality we can assume that $a$ satisfies condition~\Ref{eq_cond_old_multilin}. In particular, any multilinear path in $\supp(a)$ satisfies condition~\Ref{eq_cond_old_multilin}.

Assume that conditions~a), b), c) are not valid. Applying Corollary~\ref{cor_trAB0}, we can see that up to permutation of indices of $b_1,b_2,b_3$ one of the following possibilities holds:
\begin{enumerate}
\item[1)]  $\#(\cap_{i=1}^3 \Ver{b_i})\geq1$ and $\#(\Ver{b_1}\cap\Ver{b_2})\geq2$; 

\item[2)]  $\cap_{i=1}^3 \Ver{b_i}=\emptyset$ and $\Ver{b_i}\cap\Ver{b_j}$ is not empty for all $i,j$; 

\item[3)]  $\Ver{b_1}\cap\Ver{b_2}=\{u_1,u_2\}$, $\Ver{b_2}\cap\Ver{b_3}=\{v_1,v_2\}$,  $\Ver{b_1}\cap\Ver{b_3}=\emptyset$ for pairwise different vertices $u_1,u_2,v_1,v_2\in\Q_0$; moreover, $b_2\sim c_1 c_2$ for paths $c_1,c_2$ with $c_i',c_i''\in\{u_1,u_2\}$ and $\Ver{c_i}\cap\Ver{b_3}\neq\emptyset$ for $i=1,2$. 
\end{enumerate} 
We claim that there is a multilinear path $e$ in $\supp(a)$ such that $\tr(e)\equiv0$. By Lemma~\ref{lemma_multilinear_path}, this claim implies $\tr(a)\equiv0$; a contradiction.
To prove the claim, we consider the above mentioned cases.  

Case 1). We have $b_1\sim c_1 c_2$ and $b_2\sim d_1 d_2$ for paths $c_1,c_2,d_1,d_2$ in $\Q^{\ast}$ such that we have the following picture in $\Q^{\ast}$:
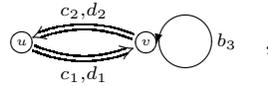
$$\xymatrix@C=1.3cm@R=1.3cm{ %
\vtx{u}\ar@2@/^/@{<-}[r]^{c_2,d_2} \ar@2@/_/@{->}[r]_{c_1,d_1} & \vtx{v}\\}%
\loopL{0}{0}{b_3}
\qquad,$$%

\noindent{}where $u\neq v$. By relation~($R_2$) from Lemma~\ref{lemma_two_rel}, $\tr(e)\equiv0$ for the multilinear path $e=b_3 c_1 c_2 d_1 d_2$. 

Case 2). We have $b_i\sim b_{i1} b_{i2}$ for paths $b_{i1},b_{i2}$ in $\Q^{\ast}$ ($1\leq i\leq 3$) such that we have the following picture in $\Q^{\ast}$:
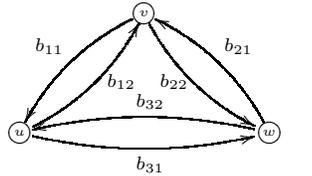
$$ 
\xymatrix@C=1.3cm@R=1.3cm{ %
&\vtx{v} %
\ar@/_/@{->}[ld]_{b_{11}} \ar@/^/@{<-}[ld]
\ar@/^/@{<-}[rd]^{b_{21}} \ar@/_/@{->}[rd] & \\
\vtx{u} \ar@/_/@{->}[rr]_{\;\;b_{31}} \ar@/^/@{<-}[rr]^{\;\;b_{32}} &&\vtx{w}  }%
\begin{picture}(0,0)(0,0)
\put(-67,-25){$\scriptstyle b_{12}$} 
\put(-47,-25){$\scriptstyle b_{22}$} 
\put(10,-25){.} 
\end{picture}
$$

\noindent{}where $u, v, w$ are pairwise different. By Lemma~\ref{lemma_trAB0}, $\tr(e)\equiv0$ for the multilinear path $e=b_{11} b_{21} b_{31} b_{32} b_{22} b_{12}$.  

Case 3). We have $b_1\sim b_{11} b_{12}$, $b_2\sim d_1\cdots d_4$, and $b_3\sim b_{31} b_{32}$ for paths $b_{1i},b_{3i},d_j$ in $\Q^{\ast}$ ($i=1,2$, $1\leq j\leq 4$) such that up to permutations of vertices $v_1,v_2$ and $u_1,u_2$ we have the following picture:
$$ 
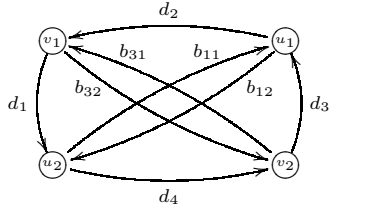
\xymatrix@C=2.6cm@R=1.3cm{ %
\vtx{v_1} \ar@/^/@{<-}[r]^{d_2} \ar@/^/@{<-}[rd] \ar@/_/@{->}[rd] \ar@/_/@{->}[d]_{d_1} &
\vtx{u_1} \ar@/^/@{<-}[d]^{d_3} \ar@/_/@{<-}[ld] \ar@/^/@{->}[ld] \\  
\vtx{u_2} \ar@/_/@{->}[r]_{d_4} & \vtx{v_2} }%
\begin{picture}(0,0)(0,0)
\put(-55,-3){$\scriptstyle b_{11}$} 
\put(-35,-17){$\scriptstyle b_{12}$} 
\put(-83,-3){$\scriptstyle b_{31}$} 
\put(-100,-17){$\scriptstyle b_{32}$} 
\put(10,-25){.} 
\end{picture}
$$

\noindent{} By Lemma~\ref{lemma_trAB0}, $\tr(e)\equiv0$ for the multilinear path %
$e=b_{11} d_4^{\ast} b_{31}^{\ast} d_2 \cdot d_3 b_{32} d_1^{\ast} b_{12}$.
\end{proof}

\begin{lemma}\label{lemma_chain of_cycles}
Let $a$ be a multilinear path in $\Q^{\ast}$. Assume that some decomposition of $a$ into primitive closed paths contains pairwise non-equivalent $b_1,\ldots,b_r$ $(r>3)$ such that $\#(\Ver{b_i}\cap\Ver{b_j})\neq\emptyset$ if and only if $|i-j|\leq 1$ or $i,j\in\{1,r\}$. Then $\tr(a)\equiv0$. 
\end{lemma}
\begin{proof} Let $\tr(a)$ be decomposable. As in the proof of Lemma~\ref{lemma_possibilities}, we can assume that $a$ satisfies condition~\Ref{eq_cond_old_multilin}. Without loss of generality we can assume that case~b) from Lemma~\ref{lemma_possibilities} holds for $b_i,b_{i+1},b_{i+2}$ for all $1\leq i\leq r$, where we set $b_{r+1}=b_1$ and $b_{r+2}=b_2$. Further we proceed as in case~2) from the proof of Lemma~\ref{lemma_possibilities}. Namely, it is not difficult to see that $b_i\sim c_i d_i$ for paths $c_i,d_i$ in $\Q^{\ast}$ such that 
\begin{enumerate}
\item[$\bullet$] $c_{i}''\in\Ver{b_i}\cap\Ver{b_{i+1}}$ for all $1\leq i\leq r$;

\item[$\bullet$] $c=c_1\cdots c_r$ and $d=d_r\cdots d_1$ are closed paths in $\Q^{\ast}$. 
\end{enumerate}
Since $cd$ is a multilinear path and $r\geq3$ we have that $\tr(cd)\equiv0$ by Lemma~\ref{lemma_trAB0}. Lemma~\ref{lemma_multilinear_path} implies a contradiction.
\end{proof}

\begin{lemma}\label{lemma_zero_type}
Let $a$ be a multilinear path in $\Q^{\ast}$ and $a$ is not admissible. Then $\tr(a)\equiv0$.
\end{lemma}
\begin{proof}
Assume that $\tr(a)$ is indecomposable.  Let $\{b_1,\ldots,b_s\}$ be some decomposition of $a$ into primitive closed paths and a diagram $D$ be the type of $a$ with respect to $\{b_1,\ldots,b_s\}$. Applying relation~($R_1$) from Lemma~\ref{lemma_two_rel}, Corollary~\ref{cor_trAB0} and Lemmas~\ref{lemma_possibilities},~\ref{lemma_chain of_cycles} to $b_1,\ldots,b_s$, we can see that $a$ is admissible with respect to $\{b_1,\ldots,b_s\}$; a contradiction.   
\end{proof}

For any $a=a_1\cdots a_s\in\path{\Q^{\ast}}$ with $a_i\in\Q_1^{\ast}$ we set $L(a)=\{a_1\}$ and $R(a)=\{a_2,\ldots,a_s\}$. Note that if $a$ is an empty path, then $L(a)=R(a)=\emptyset$; if $a\in\Q_1^{\ast}$, then $R(a)=\emptyset$.

\begin{lemma}\label{lemma_nonzero_type}
Let $a$ be an admissible multilinear path in $\Q^{\ast}$ that satisfies condition~\Ref{eq_cond_old_multilin}. Then $\tr(a)$ is indecomposable.
\end{lemma}
\begin{proof} 
Let $\{b_1,\ldots,b_s\}$ be such decomposition of $a$ into primitive closed paths that $a$ is admissible with respect to this decomposition. Assume that a diagram $D$ is the type of $a$ with respect to $\{b_1,\ldots,b_s\}$. Denote by $v_i$ the vertex of $D$ corresponding to $b_i$ (see Definition~\ref{def_type}).

We prove the lemma by induction on $s\geq1$. Let $s=1$. Assume that $\tr(a)\equiv0$. Then $\tr(a)=\sum_i \al_i \prod_j \tr(c_{ij})$, where $\al_i\in\FF$ and $c_{ij}$ is a closed path in $\Q^{\ast}$ with $\sum_j\mdeg(c_{ij})=\mdeg(a)$ and $\deg(c_{ij})<\deg(a)$. Since $\tr(a)\neq0$, we have $c_{ij}\sim a$; a contradiction.    

Let $s>1$. We shrink all triangles of $D$ and obtain a tree (see Remark~\ref{remark_shrink}). Considering  all leafs of this tree, we can see that one of the following cases holds.
\begin{enumerate}
\item[1.] There are $1\leq i,j\leq s$ ($i\neq j$) such that $v_i$ and $v_j$ are connected by means of an edge marked with $1$ and $\Ver{b_i}\cap\Ver{b_q}\neq\emptyset$ if and only if $q=i$ or $q=j$. Hence for some $c\sim b_i$ and $x,y\in\Arr{b_j}$ we have
$$
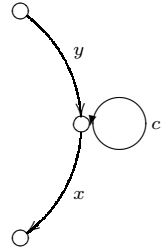
\xymatrix@C=0.6cm@R=1.3cm{ %
\vtx{}\ar@/^/@{->}[rd]^{y} & \\
& \vtx{} \\
\vtx{}\ar@/_/@{<-}[ru]_{x} & \\}
\loopL{0}{43}{c}
$$

\item[2.] There are pairwise different $1\leq i,j,k\leq s$ such that $v_i,v_j,v_k$ are vertices of a triangle and for $p=i,j$ we have $\Ver{b_p}\cap\Ver{b_q}\neq\emptyset$ if and only if $q\in\{i,j,k\}$. Hence for some $c\sim b_i$, $d\sim b_j$ and $x,y\in\Arr{b_k}$ we have
$$
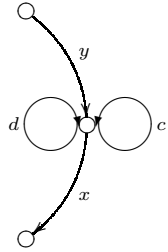
\xymatrix@C=0.6cm@R=1.3cm{ %
\vtx{}\ar@/^/@{->}[rd]^{y} & \\
& \vtx{} \\
\vtx{}\ar@/_/@{<-}[ru]_{x} & \\}
\loopL{0}{43}{c}
\loopR{48}{43}{d}
$$

\item[3.] There are $1\leq i,j\leq s$ ($i\neq j$) such that $v_i,v_j$ are connected by means of an edge marked with $2$ and $\Ver{b_i}\cap\Ver{b_q}\neq\emptyset$ if and only if $q=i,j$. Hence for some $c=c_1 c_2\sim b_i$, arrows $x,y\in\Arr{b_j}$, and a path $d$ in $\supp(b_j)$ we have
$$
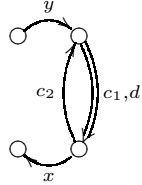
\xymatrix@C=0.6cm@R=1.3cm{ %
\vtx{}\ar@/^/@{->}[r]^{y} & \vtx{} \ar@2@/^/@{->}[d]^{c_1,d} \ar@/_/@{<-}[d]_{c_2} \\
\vtx{}\ar@/_/@{<-}[r]_{x} &  \vtx{}\\}
$$
\end{enumerate} 

Denote  $I_2=\left(
\begin{array}{cc}
1& 0 \\
0& -1\\
\end{array}
\right)$. Define a homomorphism $\Phi: SI(\Q)\to SI(\Q)$ of $\FF$-algebras as follows (see Remark~\ref{remark_Phi}): for every $z\in\Q_1^{\ast}$ we set
\begin{enumerate}%
\item[$\bullet$] in case~1 we have $\Phi(X_z)=\left\{
\begin{array}{rl}
I_2,& \text{if } z\in L(c)\\
E,& \text{if } z\in R(c)\\
X_z I_2,& \text{if } z=x \\
X_z,& \text{otherwise }\\
\end{array}
\right.$;

\item[$\bullet$] in case~2 we have $\Phi(X_z)=\left\{
\begin{array}{rl}
I_2,& \text{if } z\in L(c)\\
J_2,& \text{if } z\in L(d)\\
E,& \text{if } z\in R(c)\cup R(d)\\
X_z I_2 J_2,& \text{if } z=x \\
X_z,& \text{otherwise }\\
\end{array}
\right.$;

\item[$\bullet$] in case~3 we have $\Phi(X_z)=\left\{
\begin{array}{rl}
I_2,& \text{if } z\in L(c_1)\\
J_2,& \text{if } z\in L(c_2)\\
E,& \text{if } z\in R(c_1)\cup R(c_2)\cup\Arr{d}\\
X_z X_d I_2 J_2,& \text{if } z=x \\
X_z,& \text{otherwise }\\
\end{array}
\right..$
\end{enumerate}

If $\tr(a)\equiv 0$, then $\tr(a)=\sum_q \al_q f_q$ for $\al_q\in\FF$ and some products $f_q$ of at least two traces. Note that $\tr(I_2)=\tr(J_2)=\tr(I_2 J_2)=0$, $I_2^2=E$, $J_2^2=-E$, and $I_2 J_2 = -J_2 I_2$. Thus, there is a multilinear path $e$ in $\Q^{\ast}$  such that $\Phi(\tr(a))\equiv\pm \tr(e)$ and $\mdeg(e)=\mdeg(a)- \mdeg(c)$ in cases~1 and~3 and $\mdeg(e)=\mdeg(a)- \mdeg(c) - \mdeg(d)$ in case~2. On the other hand, $\Phi(f_q)$ is either zero or a product of at least two traces of closed paths. Therefore,
$\tr(e)\equiv0$. Since there is a decomposition of $e$ into primitive closed paths that consists of $s-1$ or $s-2$ paths and $e$ is admissible with respect to the mentioned decomposition, induction hypothesis implies a contradiction.
\end{proof}

Denote by $D(\Q)$ the set of all maps $\de:\Q^{\ast}_1\to\{0,1\}$ such that if $\de(x)=1$, then $x$ is a loop.  Given $\de\in D(\Q)$, we define a homomorphism $\Psi_{\de}: SI(\Q)\to \FF[H(\Q,(2,\ldots,2))]$ of $\FF$-algebras as follows: for every $z\in\Q_1^{\ast}$ we set 
$$\Psi_{\de}(X_z)=\left\{
\begin{array}{rl}
X_z,& \text{if } \de(z)=0 \\
X_z - \tr(X_z) E_{22},& \text{if } \de(z)=1 \\
\end{array}
\right.,$$%
where $E_{22}=\left(
\begin{array}{cc}
0 & 0 \\
0 & 1\\
\end{array}
\right)$ (see Remark~\ref{remark_Phi} for details). If particular, for $z\in\Q_1$ with $\de(z)=1$ we have 
$\Psi_{\de}(X_z)=\left(
\begin{array}{cc}
x^z_{11} & x^z_{12} \\
x^z_{21} & -x^z_{11}\\
\end{array}
\right)$, %
where $x_{ij}^z\in\FF[H(\Q,(2,\ldots,2))]$ (see Section~\ref{section_intro}). As above, for $f\in SI(\Q)$ we say that $\Psi_{\de}(f)$ is {\it decomposable} and write $\Psi_{\de}(f)\equiv0$ if $\Psi_{\de}(f)$ is a polynomial in elements from  $\Psi_{\de}(SI(\Q))$ of strictly less degree or $\Psi_{\de}(f)\in\FF$.

We say that $s$ is a {\it complexity} of a multilinear path $a$ in $\Q^{\ast}$ if there is a decomposition $\{b_1,\ldots,b_s\}$ of $a$ into primitive closed paths. Note that $a$ can have several pairwise different complexities. 

\remark{}\label{remark_double arrow} 
If $a\in\SetI$ and $b$ is a closed primitive path in $\supp(a)$ and 
$\deg_{x}(b)+\deg_{x^{\ast}}(b)>0$ for an $x\in\Q_1$, then $b\sim x x^{\ast}$ and $b$ is called a {\it double arrow} of $a$.  
\bigskip{}
 
The next lemma is a generalization of Lemma~\ref{lemma_nonzero_type}.

\begin{lemma}\label{lemma_S1}
If $a\in \SetI$ is not a loop and $\de\in D(\Q)$, then $\Psi_{\de}(\tr(a))$ is indecomposable. In particular, $\tr(a)$ is indecomposable. 
\end{lemma}
\begin{proof}  We prove the lemma by induction on complexity of $a$. For short, we write $f^{\Psi}$ for $\Psi_{\de}(f)$, where $f\in SI(\Q)$. 

If one is a complexity of $a$, then $\de(x)=0$ for all $x\in\Arr{a}$. We obtain  the required statement in the same way as in the proof of Lemma~\ref{lemma_nonzero_type}.  

We assume that $b_1, b_2, b_3\in\Q^{\ast}$ are loops in a $v\in\Q_0$ and $\de(b_i)=1$ for $1\leq i\leq 3$. We claim that 
\begin{eq}\label{eq_claim}
\tr(b_1b_2)^{\Psi}\not\equiv0 \text{ and } \tr(b_1b_2b_3)^{\Psi}\not\equiv0.
\end{eq}%
Let $\tr(b_1 b_2 b_3)^{\Psi}\equiv0$. Since $\tr(b_i)^{\Psi}=0$, we obtain $\tr(b_1 b_2 b_3)^{\Psi}=0$. But the last equality is not valid; a contradiction. In the same way we can see that $\tr(b_1b_2)^{\Psi}$ is indecomposable.

Let $\{b_1,\ldots,b_s\}$ be some decomposition of $a$ into primitive closed paths in $\Q^{\ast}$ and $s>1$. Then we can see that case $1$, $2$ or $3$ from the proof of Lemma~\ref{lemma_nonzero_type} holds. In what follows, we use notations from the proof of Lemma~\ref{lemma_nonzero_type}. By Definition~\ref{def_tree_path}, we have that
\begin{enumerate}
\item[$\bullet$] $c$ is not a double arrow in case~1, 

\item[$\bullet$] $c$ and $d$ are not double arrows in case~2,

\item[$\bullet$] $c$ and $b_j$ are not double arrows in case~3.
\end{enumerate}
We set $b=b_j$ in cases~1,3 and $b=b_k$ in case~2. 

Let $b$ be not a double arrow. If $b$ is not a loop, then we define $\Phi$ in the same way as in the proof of Lemma~\ref{lemma_nonzero_type} and apply induction hypothesis to complete the proof. If $b$ is a loop, then either case~1 or~2 holds and $\tr(a)^{\Psi}=\tr(bc)^{\Psi}$ in case~1 and $\tr(a)^{\Psi}\equiv\pm\tr(bcd)^{\Psi}$ in case~2, where $b,c$ are closed paths with $\Arr{b}\cap\Arr{c}=\emptyset$ in case~1 and similarly in case~2. Obviously,~\Ref{eq_claim} implies $\tr(a)^{\Psi}\not\equiv0$.

Let $b$ be a double arrow. Therefore,  $b\sim x x^{\ast}$ and $y=x^{\ast}$. Without loss of generality, we can assume that $x\in\Q_1$.

Case~1. Define a homomorphism $\Phi: SI(\Q)\to \FF[H(\Q,(2,\ldots,2))]$ of $\FF$-algebras as follows: for every $z\in\Q_1^{\ast}$ we set
$$\Phi(X_z)=\left\{
\begin{array}{rl}
E,& \text{if } z=x \text{ or }z=x^{\ast}\\
X_z,& \text{otherwise }\\
\end{array}
\right..$$ 
We remove the arrow $x$ from $\Q$ and glue vertices $x'$ and $x''$ together. Denote the resulting quiver by $\Gamma$. We also remove arrows $x,x^{\ast}$ from $a$ and obtain a new path $e$ in $\Gamma^{\ast}$ satisfying $\tr(e)=\Phi(\tr(a))$. Moreover, $e$ is an admissible tree path in $\Gamma^{\ast}$ and $e$ is not a loop. Let $\tr(a)^{\Psi}\equiv0$. Applying relation (C), we obtain $\tr(e)^{\Psi}\equiv0$. Since a complexity of $e$ is equal to $s-1$, induction hypothesis implies a contradiction.

Case~2. For $X,C,D\in\FF^{2\times 2}$ we define a homomorphism $\Phi=\Phi_{X,C,D}: SI(\Q)\to \FF[H(\Q,(2,\ldots,2))]$ of $\FF$-algebras as follows: for every $z\in\Q_1^{\ast}$ we set
$$\Phi(X_z)=\left\{
\begin{array}{rl}
X,& \text{if } z=x \\
X^{\ast},& \text{if } z=x^{\ast} \\
C,& \text{if } z=L(c) \\
D,& \text{if } z=L(d) \\
E,& \text{if } z\in R(c)\cup R(d) \\
X_z,& \text{otherwise }\\
\end{array}
\right..$$ 
We remove the arrows $\{x,x^{\ast}\}\cup\Arr{c}\cup\Arr{d}$ from $\Q^{\ast}$, add a new loop $y$ in the vertex $x''$, and glue vertices $x'$ and $x''$ together. As the result of this procedure, we obtain a quiver $\Gamma^{\ast}$ for some quiver $\Gamma$. We remove arrows $\{x^{\ast}\}\cup\Arr{c}\cup\Arr{d}$ from $a$ and substitute $y$ for $x$. As the result, we obtain a path $e$ in $\Gamma^{\ast}$. Note that $e$ is an admissible tree path in $\Gamma^{\ast}$ and $e$ is not a loop. We set $\de(y)=1$. Thus, we can consider $\de$ as a map $\Gamma^{\ast}_1\to \{0,1\}$. 

Let $\tr(a)^{\Psi}\equiv0$. Applying relations~$(\rm A')$ and~$(\rm C)$, we obtain
$$\tr(a)^{\Psi}=\sum\nolimits_i \al_i \tr(xcdx^{\ast}z_i)^{\Psi} f_i + \sum\nolimits_j \be_j h_j,$$
where $\al_i,\be_j\in\FF$, $f_i,h_j\in \Psi_{\de}(SI(\Q))$, $z_i$ is a closed path in $\Q^{\ast}$, and $h_j$ does not contain neither $\tr(xcdx^{\ast}z)^{\Psi}$ nor $\tr(xdcx^{\ast}z)^{\Psi}$ as a factor for any closed path $z$ in $\Q^{\ast}$. Assume that
$$\tr(C)=\tr(D)=\tr(CD)=0 \text{ and } \tr(X C D X^{\ast})=0.$$
Since $\Phi(\sum_j \be_j h_j)=0$ and $\tr(Y)=0$ for $Y=XCDX^{\ast}$, we obtain that 
$\tr(e)^{\Psi}|_{X_y\to Y} = \sum\nolimits_i \al_i \tr(Y X_{z_i})^{\Psi} f_i$. Lemma~\ref{lemma_null_sled} (see below) implies that $\tr(e)^{\Psi} = \sum\nolimits_i \al_i \tr(X_y X_{z_i})^{\Psi} f_i$. Thus, $\tr(e)^{\Psi}\equiv0$.  Since a complexity of $e$ is equal to $s-2$, induction hypothesis implies a contradiction. 
\end{proof}

\begin{lemma}\label{lemma_null_sled}
Assume that $Y\in\FF^{2\times 2}$ satisfies $\tr(Y)=0$. Then there are $X,C,D\in\FF^{2\times 2}$ such that $\tr(C)=\tr(D)=\tr(CD)=0$ and $Y=XCDX^{\ast}$.
\end{lemma}
\begin{proof}
We set $Y=\left(
\begin{array}{cc}
y_1& y_2 \\
y_3& -y_1\\
\end{array}
\right)$.

Let $y_2$ and $y_3$ be non-zero or $y_2=y_3=0$. We take $X=E$,  
$C=\left(
\begin{array}{cc}
-c_1& c_2 c_3 \\
-c_3& c_1\\
\end{array}
\right)$ 
and
$D=\left(
\begin{array}{cc}
0& c_2 \\
1& 0\\
\end{array}
\right)$ for $c_1,c_2,c_3\in\FF$.
If $y_2=y_3=0$, then we consider $c_1=0$, $c_2=1$, $c_3=y_1$ and obtain the required. If $y_2$ and $y_3$ are non-zero, then we consider $c_1=y_3$, $c_2=-y_2/y_3$, $c_3=-y_1 y_3/y_2$ and the required statement follows. 
  
Let $y_2=0$ and $y_3\neq0$. Then matrices 
$X=\left(
\begin{array}{cc}
1& y_1/y_3 \\
0& 1\\
\end{array}
\right)$,
$C=\left(
\begin{array}{cc}
-y_3& 0 \\
0& y_3\\
\end{array}
\right)$, and
$D=\left(
\begin{array}{cc}
0 & -y_1^2/y_3^2 \\
1 & 0\\
\end{array}
\right)$ satisfy the required property.

Let $y_2\neq0$ and $y_3=0$. Then matrices 
$X=\left(
\begin{array}{cc}
1& 0 \\
-y_1/y_2& 1\\
\end{array}
\right)$,
$C=\left(
\begin{array}{cc}
0& 1 \\
-y_1^2/y_2^2 &0\\
\end{array}
\right)$, and
$D=\left(
\begin{array}{cc}
-y_2 & 0 \\
0 & y_2\\
\end{array}
\right)$ satisfy the required property.
\end{proof}

\begin{proof_of}{of Theorem~\ref{theo_mgs_not2}.}
Denote by $P$ the set from the formulation of the theorem. Theorem~\ref{theo_Fedotov} together with Lemmas~\ref{lemma_mdeg_tr},~\ref{lemma_not_tree_path},~\ref{lemma_zero_type} and relations~(B),~(E) show that $P$ generates $SI(\Q)$. Lemma~\ref{lemma_S1} together with the fact that all elements of $P$ have pairwise different multidegrees imply that $P$ is a minimal generating set for $SI(\Q)$.
\end{proof_of}

\section{Corollaries}\label{section_corollaries}

In this section we collect some corollaries concerning $SI(\Q)$. Let us recall that if we consider arrows of a quiver $\Q$ as an undirected edges, then $\Q$ turns into the {\it underlying graph} of $\Q$.

\begin{cor}\label{cor_orientation}
Let $\Gamma$ and $\Lambda$ be quivers with isomorphic underlying graphs. Then $SI(\Gamma)\simeq SI(\Lambda)$. Moreover, the given isomorphism preserves multidegrees.
\end{cor}
\begin{proof} Since the underlying graphs of $\Gamma$ and $\Lambda$ are isomorphic, there are isomorphisms $\varphi_0:\Gamma_0\to \Lambda_0$ and $\varphi_1:\Gamma_1\to \Lambda_1$ such that for every $a\in\Gamma_1$ we have $\{\varphi_0(a'),\varphi_0(a'')\}=\{\varphi_1(a)',\varphi_1(a)''\}$. 

Define the map $\varphi:\Gamma^{\ast}_1\to \Lambda^{\ast}_1$ as follows: 
$$\varphi(a)=\left\{\begin{array}{cl}
\varphi_1(a),& \text{if } \varphi_0(a')=\varphi_1(a)'\\
\varphi_1(a)^{\ast},&\text{otherwise}\\
\end{array}\right.\qquad
\text{and }\qquad \varphi(a^{\ast})=\varphi(a)^{\ast},$$
where $a\in\Gamma_1$. Given a path $a=a_1\cdots a_s$ in $\Gamma^{\ast}$, where $a_1,\ldots,a_s\in\Gamma_1^{\ast}$, we write $\varphi(a)$ for $\varphi(a_1)\cdots\varphi(a_s)$. Obviously, if $a$ is a closed path in $\Gamma^{\ast}$, then $\varphi(a)$ is a closed path in $\Lambda^{\ast}$.

Define the homomorphism of algebras $\Phi:\FF[H(\Gamma,(2,\ldots,2))]\to \FF[H(\Lambda,(2,\ldots,2))]$ as follows: $\Phi(x_{ij}^a)$ is $(i,j)^{\rm th}$ entry of $X_{\varphi(a)}$, where $a\in\Gamma_1$. Since $\Phi(\det(X_a))=\det(X_{\varphi(a)})$ and  $\Phi(\tr(X_b))=\tr(X_{\varphi(b)})$ for an $a\in\Gamma_1$ and a closed path $b\in\Gamma^{\ast}$, Theorem~\ref{theo_Fedotov} implies that the restriction of $\Phi$ to $SI(\Gamma)$ is an epimorphism $SI(\Gamma)\to SI(\Lambda)$. Considering $\varphi_0^{-1}$ instead of $\varphi_0$ and $\varphi_1^{-1}$ instead of $\varphi_1$ and repeating the above reasoning, we construct an epimorphism $\Psi:SI(\Lambda)\to SI(\Gamma)$ such that $\Phi\circ \Psi$ is the identity map. 
\end{proof}

\begin{remark} Corollary~\ref{cor_orientation} does not hold for an arbitrary dimension vector $\n$. As an example, assume that 
$$\Gamma:\quad
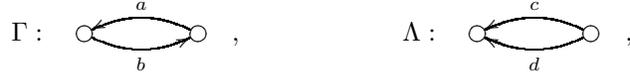
\xymatrix@C=1.3cm@R=1.3cm{ %
\vtx{} \ar@/^/@{<-}[r]^{a} \ar@/_/@{->}[r]_{b} & \vtx{}    \\
}\quad,\qquad\qquad\qquad%
\Lambda:\quad
\xymatrix@C=1.3cm@R=1.3cm{ %
\vtx{} \ar@/^/@{<-}[r]^{c} \ar@/_/@{<-}[r]_{d} & \vtx{}    \\
}\quad,
$$
and $\n=(3,3)$. Then $f=\tr(X_a X_b)\in SI(\Gamma,\n)$ and $\deg(f)=2$. On the other hand,  for any $h\in SI(\Lambda,\n)$ with $h\not\in\FF$ we have $\deg(h)\geq3$. Therefore, there is no an isomorphism between $SI(\Gamma,\n)$ and $SI(\Lambda,\n)$ that preserves multidegrees, but the underlying graphs of $\Gamma$ and $\Lambda$ are isomorphic. 
\end{remark}

\begin{cor}\label{cor_MGS}
Let $\{f_1,\ldots,f_r\}$ and $\{h_1,\ldots,h_s\}$ be minimal $\NN^{\#\Q_1}$-homogeneous generating sets for $SI(\Q)$.
Then $r=s$ and there is a $\pi\in S_r$ and non-zero $\al_1,\ldots,\al_r\in\FF$ such that $f_1\equiv \al_1 h_{\pi(1)},\ldots,f_r\equiv \al_r h_{\pi(r)}$.
\end{cor}
\begin{proof} Theorem~\ref{theo_Fedotov} and Lemma~\ref{lemma_mdeg_tr} together with relation~(E) of Lemma~\ref{lemma_relations} imply that the dimension of every $\NN^{\#\Q_1}$-homogeneous component of $\ov{SI(\Q)}={SI(\Q)}/{(SI(\Q)^{+})^2}$ is either $0$ or $1$. Since an $\NN^{\#\Q_1}$-homogeneous set $\{f_i\} \subseteq SI(\Q)$ is a minimal set of generators of $SI(\Q)$ if and only if $\{\ov{f_i}\}$ is a basis of $\ov{SI(\Q)}$, the proof is completed. 
\end{proof}

\begin{cor}\label{cor_relations}
Any relation $\sum_i \al_i f_i\equiv0$, where $\al_i\in\FF$ and $f_i\in SI(\Q)$ is indecomposable, is a linear combination of relations from Lemma~\ref{lemma_relations}.
\end{cor} 
\begin{proof}
In the proofs of Theorems~\ref{theo_mgs_2} and~\ref{theo_mgs_not2} we show that using relations from Lemma~\ref{lemma_relations} we can represent any semi-invariant as a linear combination of elements from the minimal generating set modulo decomposable semi-invariants. Thus, the required statement holds.
\end{proof}

\section{Examples}\label{section_example}

In this section we apply our main result to tree-like quivers and its partial cases such as tree quivers and quivers with two vertices.

\subsection{Tree-like quivers}

We say that a quiver $\Q$ is a {\it tree} if its underlying graph is a tree. A quiver $\Q$ is called a {\it tree-like} quiver if the degree of every primitive closed path in the underlying graph of $\Q$ is less than three. 

Given a quiver $\Q$ and its underlying graph $\Gamma$, denote by $\widehat{\Q}$ the graph that we obtain from $\Gamma$ as the result of the following procedure:
\begin{enumerate}
\item[$\bullet$] remove all loops from $\Gamma$;

\item[$\bullet$] for every $u,v\in\Gamma_0$ with $u\neq v$ consider the set of edges connecting $u,v$ and remove all edges from it but one.
\end{enumerate}%
Obviously, $\Q$ is a tree-like quiver if and only if $\widehat{\Q}$ is a tree. 

\example\label{ex_treelike} The following quiver $\Q$ is a tree-like quiver, since $\widehat{\Q}$ is a tree:
$$
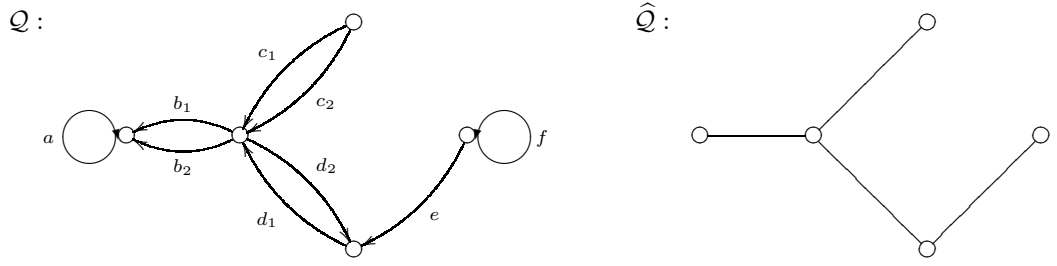
\Q: \;\; \loopR{1}{44}{a\,}
\xymatrix@C=1.3cm@R=1.3cm{ %
&& \vtx{}  \\ 
\vtx{} \ar@/^/@{<-}[r]^{b_1} \ar@/_/@{<-}[r]_{b_2} &\vtx{} 
\ar@/^/@{->}[dr]^{d_2} \ar@/_/@{<-}[dr]_{d_1} 
\ar@/^/@{<-}[ur]^{c_1} \ar@/_/@{<-}[ur]_{c_2} & &\vtx{} \ar@/^/@{->}[dl]^{e} \\
&& \vtx{} &\\
}%
\loopL{-1}{44}{f}
\qquad\qquad \widehat{\Q}:\;\; 
\xymatrix@C=1.3cm@R=1.3cm{ %
&& \vtx{}  \\ 
\vtx{} \ar@{-}[r]  &\vtx{} \ar@{-}[dr] \ar@{-}[ur] &&%
\vtx{} \ar@{-}[dl] \\
&& \vtx{} &\\
}%
$$
\bigskip

Two arrows $a,b$ in $\Q$ are {\it parallel} if $\{a',a''\}=\{b',b''\}$. Similarly, parallel edges are defined for a graph.

We say that $\theta=(\theta_v,\theta_x\,|\,v\in \widehat{\Q}_0,\,x\in \widehat{\Q}_1)$ is a {\it coloring} of $\widehat{\Q}$ if 
\begin{enumerate}
\item[$\bullet$] $\theta_v$ is a subset of loops of $\Q$ in the vertex $v$;  

\item[$\bullet$] $\theta_x$ is a subset of arrows of $\Q$ that are parallel to $x$. 
\end{enumerate}  
Given $a\in\Q_1$, we write $a\in\theta$ if $a\in\theta_v$ or $a\in\theta_x$ for some $v$ and $x$. If $a\not\in\theta$ for all $a\in\Q_1$, then $\theta$ is called empty. Let us remove from $\widehat{\Q}$ 
\begin{enumerate}
\item[$\bullet$] all edges $x$ with empty $\theta_x$,

\item[$\bullet$] all vertices $v$ with empty set $\{y\in\widehat{Q}_1\,|\, v \text{ is a vertex of } y \text{ and } \theta_y \text{ is not empty}\}$ 
\end{enumerate} 
and denote the resulting graph by $\widehat{\Q}(\theta)$. If $\widehat{\Q}(\theta)$ is connected, then $\theta$ is called {\it connected}. 

Let $\Q$ be a tree-like quiver. A coloring $\theta$ is called {\it good} if 
\begin{enumerate}
\item[a)] $\theta$ is connected and $\theta$ is not empty;

\item[b)] for every $x\in \widehat{\Q}_1$ we have that either $\#\theta_x$ is even or $\#\theta_x=1$;

\item[c)] if $\#\theta_x=1$ for an $x\in \widehat{\Q}_1$ and one of two vertices of $x$ is a leaf $v$ of the tree $\widehat{\Q}(\theta)$, then  $\theta_v$ is not empty.  
\end{enumerate}
It is well-known that if $\Gamma$ is a connected graph such that for all $u,v\in\Gamma_0$ with $u\neq v$ number of edges connecting $u$ and $v$ is even, then there is a closed path in $\Gamma$ containing every edge of $\Gamma$ one time exactly. Hence we obtain that for any good coloring $\theta$ there is a closed path $b_{\theta}$ in $\Q^{\ast}$ such that for every $a\in \Q_1$ we have
$$\deg_{a}(b_{\theta})+\deg_{a^{\ast}}(b_{\theta})=
\left\{
\begin{array}{rl}
0,& \text{if }  a\not\in\theta\\
1,& \text{if }  a\in\theta \text{ and } a \text{ is a loop}\\
1,& \text{if }  a\in\theta,\; a \text{ is not a loop and }\#\theta_x>1\\
2,& \text{if }  a\in\theta,\; a \text{ is not a loop and }\#\theta_x=1\\
\end{array}
\right.,$$
where in the $3^{\rm rd}$ and $4^{\rm th}$ cases $x$ stands for the only edge in $\widehat{\Q}$ parallel to $a$. Note that 
$$\deg_{a}(b_{\theta})=\deg_{a^{\ast}}(b_{\theta})=1$$
in the last case. 

\example\label{ex_coloring} Let $\Q$ be the tree-like quiver from Example~\ref{ex_treelike}. Then the following coloring $\theta$ of $\widehat{\Q}$ is good, where we write down $\theta_x$ ($\theta_v$, respectively) near arrow $x$ (vertex $v$, respectively) of $\widehat{\Q}$. We also depict $\widehat{\Q}(\theta)$:
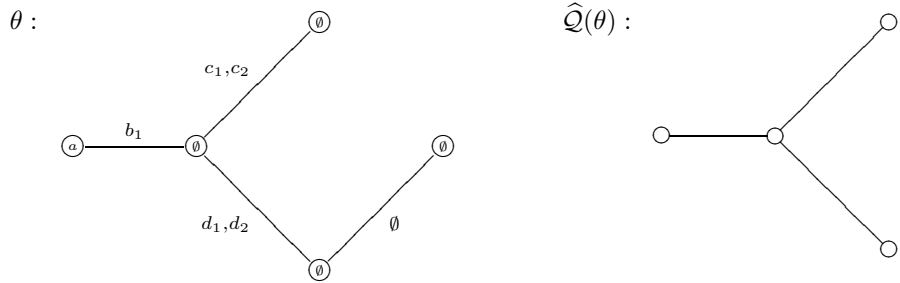
$$
\theta:\;\; 
\xymatrix@C=1.3cm@R=1.3cm{ %
&& \vtx{\emptyset}  \\ 
\vtx{a} \ar@{-}[r]^{b_1}  &\vtx{\emptyset} \ar@{-}[dr]_{d_1,d_2} \ar@{-}[ur]^{c_1,c_2} &&%
\vtx{\emptyset} \ar@{-}[dl]^{\emptyset} \\
&& \vtx{\emptyset} &\\
}%
\qquad\qquad\widehat{\Q}(\theta):\;\; 
\xymatrix@C=1.3cm@R=1.3cm{ %
&& \vtx{}  \\ 
\vtx{} \ar@{-}[r]  &\vtx{} \ar@{-}[dr] \ar@{-}[ur] && \\
&& \vtx{} &\\
}%
$$
We can assume that $b_{\theta}$ is the following closed path in $\Q^{\ast}$: $a b_1 c_1 c_2^{\ast} d_1 d_2 b_1^{\ast}$.

\begin{lemma}\label{lemma_tree_like_2}
Assume that $\Char{\FF}=2$ and $\Q$ is a tree-like quiver. Then the following set is a minimal generating set for $SI(\Q)$: 
$$\det(X_a),\,\tr(X_{b_\theta}),$$%
where $a$ ranges over $\Q_1$ and $\theta$ ranges over good colorings of $\widehat{\Q}$.  
\end{lemma}
\begin{proof} Since $\theta$ is a good coloring, $b_{\theta}$ is a tree path in $\Q^{\ast}$ (see Definition~\ref{def_tree_path}).  Obviously, if $\theta$ and $\vartheta$ are pairwise different colorings of $\widehat{\Q}$, then multidegrees of $b_{\theta}$ and $b_{\vartheta}$ are different.  Theorem~\ref{theo_mgs_2} completes the proof.
\end{proof}
 
\begin{lemma}\label{lemma_tree_like_not2}
Assume that $\Char{\FF}\neq 2$ and $\Q$ is a tree-like quiver. Then the following set is a minimal generating set for $SI(\Q)$: 
$$\det(X_a),\,\tr(X_{b_\theta}),$$%
where $a$ ranges over $\Q_1$ and $\theta$ ranges over good colorings of $\widehat{\Q}$ satisfying one the following conditions:
\begin{enumerate}
\item[a)] there is an $x\in\widehat{\Q}_1$ such that $\#\theta_x=4$ and $\theta_v,\theta_y$ are empty for all $v\in \widehat{Q}_0$ and $y\in\widehat{Q}_1$ with $y\neq x$; 

\item[b)] $\#\theta_x\leq2$ for all $x\in\widehat{Q}_1$ and for every $v\in\widehat{\Q}_0$ we have
$$\#\theta_v + \#\{y\in\widehat{\Q}_1\,|\, v \text{ is a vertex of } y \text{ and } \theta_y \text{ is not empty}\}\leq 3.$$
\end{enumerate}
\end{lemma}
\begin{proof} Assume that case~a) holds. Then $b_{\theta}$ is a tree path and its type 
with respect to any decomposition into primitive closed paths is the diagram 
$$\xymatrix@C=1.3cm@R=2cm{ %
\vtx{} \ar@/^/@{-}[r]^{2} & \vtx{}}\quad.%
$$%
(see Definitions~\ref{def_tree_path} and~\ref{def_type}). Definition~\ref{def_type_admissible} implies that $b_{\theta}$ is admissible. 

Assume that case~b) holds. Then $b_{\theta}$ is a tree path with unique decomposition into primitive closed paths up to $\sim$-equivalence and permutations of decomposition's elements (see Definition~\ref{def_decomposition}). The type of $b_{\theta}$ with respect to this decomposition is a tree $D$ such that all edges of $D$ are marked with $1$. It is not difficult to see that $b_{\theta}$ is admissible.  

We have shown that the set from the formulation of the lemma lies in the minimal generating set for $SI(\Q)$ from Theorem~\ref{theo_mgs_not2}. Similarly we obtain the inverse inclusion of the sets. 
\end{proof}

\begin{cor}\label{cor_tree}
If $\Q$ is a tree, then $SI(\Q)=\FF[\det(X_a)\,|\,a\in\Q_1]$ is a polynomial algebra (i.e.~a free algebra over $\FF$).
\end{cor}
\begin{proof} Let $\Char{\FF}=2$ and $\theta$ be a good coloring of $\widehat{\Q}$.  Consider a leaf $v$ of $\widehat{\Q}(\theta)$ and an arrow $x\in\widehat{\Q}_1$ such that $v$ is one of two vertices of $x$. Since $\theta_v$ is empty, part~c) of the definition of good coloring implies that $\#\theta_x\neq1$. Thus $\theta_x$ is empty. Since $\theta$ is connected, $\theta$ is empty; a contradiction with part~a) of the definition of good coloring. Lemma~\ref{lemma_tree_like_2} completes the proof. 

Let $\Char{\FF}\neq 2$ and a good coloring $\theta$ of $\widehat{\Q}$ satisfy condition~a) or~b) from Lemma~\ref{lemma_tree_like_not2}. If condition~a) holds, then $\Q$ is not a tree. If condition~b) holds, then we use Lemma~\ref{lemma_tree_like_not2} instead of Lemma~\ref{lemma_tree_like_2} to prove the lemma similarly to the case of $\Char{\FF}=2$.
\end{proof}

\subsection{Quivers with two vertices}

Assume that $\Q$ is an arbitrary quiver with two vertices and its underlying graph is connected. Denote vertices of $\Q$ by $u$ and $v$. By Corollary~\ref{cor_orientation}, without loss of generality we can assume that there are no arrows from $u$ to $v$ in $\Q$.   Obviously, $\Q$ is a tree-like quiver. Denote arrows of $\Q$ by $x_1,\ldots,x_p$, $y_1,\ldots,y_q$, $z_1,\ldots,z_l$, where $x_i$ is a loop in $u$, $y_j$ is a loop in $v$, and $z_k$ goes from $v$ to $u$. Schematically, we depict $\Q$ as follows:
$$\loopR{0}{0}{x_1,\ldots,x_p} %
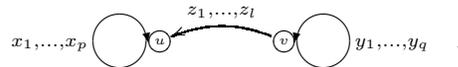
\xymatrix@C=1.3cm@R=1.3cm{ 
\vtx{u} \ar@/^/@{<-}[r]^{z_1,\ldots,z_l} & \vtx{v}\\}%
\loopL{0}{0}{y_1,\ldots,y_q}\qquad\qquad.
$$
\smallskip

\noindent{}Then the set $\SetII$ (see Section~\ref{section_results_char2}) consists of the following paths:
\begin{enumerate}
\item[a)] $x_{i_1}\cdots x_{i_r}$, where $r>0$ and $1\leq i_1<\cdots< i_r\leq p$; 

\item[b)] $y_{j_1}\cdots y_{j_s}$, where $s>0$ and $1\leq j_1<\cdots< j_s\leq q$; 

\item[c)] $x_{i_1}\cdots x_{i_r}\cdot z_k \cdot y_{j_1}\cdots y_{j_s} \cdot z_k^{\ast}$, where $r,s>0$, $1\leq k\leq l$,  $1\leq i_1<\cdots< i_r\leq p$, and $1\leq j_1<\cdots< j_s\leq q$; 

\item[d)] $x_{i_1}\cdots x_{i_r}\cdot z_{k_1} \cdot y_{j_1}\cdots y_{j_s} \cdot z_{k_2}^{\ast} \cdots z_{k_{2t-1}} z_{k_{2t}}^{\ast}$, where $r,s\geq0$, $t>0$, $1\leq k_1<\cdots< k_{2t}\leq l$, $1\leq i_1<\cdots< i_r\leq p$, and $1\leq j_1<\cdots< j_s\leq q$. \end{enumerate}  

The set $\SetI$ (see Definition~\ref{def_S1}) consists of paths from $\SetII$ that satisfy the following additional conditions, respectively:
\begin{enumerate}
\item[a)] $r\leq 3$; 

\item[b)] $s\leq 3$; 

\item[c)] $r,s\leq 2$;

\item[d)] $t\leq 2$; moreover, if $t=1$, then $r,s\leq2$; if $t=2$, then $r=s=0$.
\end{enumerate} 

Lemmas~\ref{lemma_tree_like_2} and~\ref{lemma_tree_like_not2} imply the following result.

\begin{lemma}\label{lemma_two_vertices}
A minimal generating set for $SI(\Q)$ is
\begin{enumerate}
\item[$\bullet$] $P_2=\{\det(X_a),\;\tr(X_b) \,|\, a\in\Q_1,\,b\in\SetII\}$, if 
$\Char{\FF}=2$;

\item[$\bullet$] $P_1=\{\det(X_a),\; \tr(X_b)\,|\, a\in\Q_1,\, b\in\SetI\}$, if $\Char{\FF}\neq2$.  
\end{enumerate}
\end{lemma}

\begin{remark}\label{remark_compare} Note that the generating set $P$ for $SI(\Q)$ from Theorem~\ref{theo_Fedotov} is essentially bigger than the minimal generating sets from Lemma~\ref{lemma_two_vertices}. As an example, if $p=q=l=4$, then $\#P_1=1167$, $\#P_2=2734$, but $\#P> 10^{8}$.
\end{remark}


\end{document}